\RequirePackage{fix-cm}
\documentclass[smallextended]{svjour3}       
\smartqed  
\usepackage{graphicx}
\usepackage[misc]{ifsym}
\usepackage{multirow}
\usepackage{amsmath,amssymb}
\usepackage{float}
\usepackage{url}
\usepackage[ruled]{algorithm2e}
\usepackage[top=1.0in, bottom=1.0in, left=0.9in, right=0.9in]{geometry}
\usepackage{caption}
\usepackage{algpseudocode}
\usepackage[colorlinks,linkcolor=blue,citecolor=blue]{hyperref}
\if@onethmnum
\newtheorem{assumption}[theorem]{Assumption}
\else
\newtheorem{assumption}[theorem]{Assumption}
\fi
\begin{document}
\newpage

\newpage

\title{An efficient duality-based approach for PDE-constrained sparse optimization
}


\author{Xiaoliang Song$\mathbf{^1}$  \and  Bo Chen$\mathbf{^{2}}$  \and  Bo Yu$\mathbf{^{1,3}}$}


\institute{Xiaoliang Song (\Letter)\at
           \email{songxiaoliang@mail.dlut.edu.cn}\at
           The research of this author was supported by China Scholarship Council while visiting the National University of Singapore and the National Natural Science Foundation of China (Grand No. 91230103, 11571061, 11401075)
           \and
            Bo Chen \at
           \email{chenbo@u.nus.edu}
           \and
           Bo Yu  \at
           \email{yubo@dlut.edu.cn}
           \and
           \at
           $\mathbf{^1}$ School of Mathematical Sciences,
Dalian University of Technology, Dalian, Liaoning, 116025, China
           \and
           \at
            $\mathbf{^2}$ Department of Mathematics, National University of Singapore, 10 Lower Kent Ridge Road, 119076, Singapore
           \and
           \at
           $\mathbf{^3}$ School of Mathematical and Physics Science,
Dalian University of Technology at Panjin, Panjin, Liaoning, 124200, China
}


\date{Received: date / Accepted: date}

\maketitle

\begin{abstract}
In this paper, elliptic optimal control problems involving the $L^1$-control cost ($L^1$-EOCP) is considered. To numerically discretize $L^1$-EOCP, the standard piecewise linear finite element is employed. However, different from the finite dimensional $l^1$-regularization optimization, the resulting discrete $L^1$-norm does not have a decoupled form. A common approach to overcome this difficulty is employing a nodal quadrature formula to approximately discretize the $L^1$-norm. It is clear that this technique will incur an additional error. To avoid the additional error, solving $L^1$-EOCP via its dual, which can be reformulated as a multi-block unconstrained convex composite minimization problem, is considered. Motivated by the success of the accelerated block coordinate descent (ABCD) method for solving large scale convex minimization problems in finite dimensional space, we consider extending this method to $L^1$-EOCP. Hence, an efficient inexact ABCD method is introduced for solving $L^1$-EOCP. The design of this method combines an inexact 2-block majorized ABCD and the recent advances in the inexact symmetric Gauss-Seidel (sGS) technique for solving a multi-block convex composite quadratic programming whose objective contains a nonsmooth term involving only the first block. The proposed algorithm (called sGS-imABCD) is illustrated at two numerical examples. Numerical results not only confirm the finite element error estimates, but also show that our proposed algorithm is more efficient than (a) the ihADMM (inexact heterogeneous alternating direction method of multipliers), (b) the APG (accelerated proximal gradient) method.
\end{abstract}

\section{Introduction}\label{sec:1}
In this paper, we study the following linear-quadratic elliptic PDE-constrained optimal control problem with $L^1$-control cost and piecewise box constraints on the control:
\begin{equation}\label{eqn:orginal problems}
           \qquad \left\{ \begin{aligned}
        &\min \limits_{(y,u)\in Y\times U}^{}\ \ J(y,u)=\frac{1}{2}\|y-y_d\|_{L^2(\Omega)}^{2}+\frac{\alpha}{2}\|u\|_{L^2(\Omega)}^{2}+\beta\|u\|_{L^1(\Omega)} \\
        &\qquad{\rm s.t.}\qquad Ly=u+y_r\ \ \mathrm{in}\  \Omega, \\
         &\qquad \qquad \qquad y=0\qquad\quad  \mathrm{on}\ \partial\Omega,\\
         &\qquad \qquad\qquad u\in U_{ad}=\{v(x)|a\leq v(x)\leq b,\ {\rm a.e. }\  \mathrm{on}\ \Omega\}\subseteq U,
                          \end{aligned} \right.\tag{$\mathrm{P}$}
 \end{equation}
where $Y:=H_0^1(\Omega)$, $U:=L^2(\Omega)$, $\Omega\subseteq \mathbb{R}^n$ ($n=2$ or $3$) is a convex, open and bounded domain with $C^{1,1}$- or polygonal boundary $\Gamma$; the desired state $y_d\in L^2(\Omega)$  and the source term $y_r \in L^2(\Omega)$ are given; and $a\leq0\leq b$ and $\alpha$, $\beta>0$. Moreover, the operator $L$ is a second-order linear elliptic differential operator. It is well-known that $L^1$-norm could lead to sparse optimal control, i.e. the optimal control with small support. Such an optimal control problem (\ref{eqn:orginal problems}) plays an important role for the placement of control devices \cite{Stadler}. In some cases, it is difficult or undesirable to place control devices all over the control domain and one hopes to localize controllers in small and effective regions, the $L^1$-solution gives information about the optimal location of the control devices.

Through this paper, let us suppose the elliptic PDEs involved in (\ref{eqn:orginal problems}) which are of the form
\begin{equation}\label{eqn:state equations}
\begin{aligned}
 Ly&=u+y_r \quad \mathrm{in}\  \Omega, \\
 y&=0  \qquad \quad\ \mathrm{on}\ \partial\Omega,
\end{aligned}
\end{equation}
satisfy the following assumption:
\begin{assumption}\label{equ:assumption:1}
The linear second-order differential operator $L$ is defined by
 \begin{equation}\label{operator A}
   (Ly)(x):=-\sum \limits^{n}_{i,j=1}\partial_{x_{j}}(a_{ij}(x)y_{x_{i}})+c_0(x)y(x),
 \end{equation}
where functions $a_{ij}(x), c_0(x)\in L^{\infty}(\Omega)$, $c_0\geq0$,
and it is uniformly elliptic, i.e. $a_{ij}(x)=a_{ji}(x)$ and there is a constant $\theta>0$ such that
\begin{equation}\label{equ:operator A coercivity}
  \sum \limits^{n}_{i,j=1}a_{ij}(x)\xi_i\xi_j\geq\theta\|\xi\|^2 \quad \mathrm{for\ a.a.}\ x\in \Omega\  \mathrm{and}\  \forall \xi \in \mathbb{R}^n.
\end{equation}
\end{assumption}

The weak formulation of (\ref{eqn:state equations}) is given by
\begin{equation}\label{eqn:weak form}
  \mathrm{Find}\ y\in H_0^1(\Omega):\ a(y,v)=(u+y_r,v)_{L^2(\Omega)}\quad \forall v \in H_0^1(\Omega),
\end{equation}
with the bilinear form
\begin{equation}\label{eqn:bilinear form}
  a(y,v)=\int_{\Omega}(\sum \limits^{n}_{i,j=1}a_{ji}y_{x_{i}}v_{x_{i}}+c_0yv)\mathrm{d}x,
\end{equation} or in short $ Ay=B(u+y_r)$,
where $A\in \mathcal{L}(Y,Y^*)$ is the operator induced by the bilinear form $a$, i.e., $Ay=a(y,\cdot)$ and $B\in \mathcal{L}(U,Y^*)$ is defined by $Bu=(u,\cdot)_{L^2(\Omega)}$. Since the bilinear form $a(\cdot,\cdot)$ is symmetric and $U,Y$ are Hilbert spaces, we have $A^*\in\mathcal{L}(Y,Y^*)=A$ and $B^*\in\mathcal{L}(Y,U)$ with $B^*v=v$ for any $v\in Y$.

\begin{remark}\label{more general case}
Although we assume that the Dirichlet boundary condition $y=0$ holds, it should be noted that the
assumption is not a restriction and our considerations can also carry over to the more general boundary conditions of Robin type
\begin{equation*}
  \frac{\partial y}{\partial \nu}+\gamma y=g \quad {\rm on}\  \partial\Omega,
\end{equation*}
where $g\in L^2(\partial\Omega)$ is given and $\gamma\in L^{\infty}(\partial\Omega)$ is nonnegative coefficient. Furthermore, it is assumed that the control satisfies $a\leq u \leq b$, where $a$ and $b$ have opposite signs. First, we should emphasize that this condition is required in practice, e.g., the placement of control devices. In addition, please also note, that this condition is not a restriction from the point of view of the algorithm.
If one has, e.g., $a>0$ on $\Omega$, the $L^1$-norm in $U_{ad}$ is in fact a linear
function, and thus the problem can also be handled by our method.

\end{remark}

Optimal control problems with $\alpha>0$, $\beta=0$ and their numerical realization have been studied intensively in recent papers, see e.g. \cite{Hinze,Error1,Error2,Error4,Error5,Error6} and the references cited there. Let us first comment on known results on error estimates of control constrained optimal control problems. Basic a priori error estimates were derived by Falk \cite{Error1} and Geveci \cite{Error2} where Falk considered distributed controls, while Geveci concentrated on the Neumann boundary controls. Both the authors proved optimal $L^2$-error estimates $O(h)$ for piecewise constant approximations of the control variables. Convergence results for the approximations of the controls by piecewise linear, globally continuous elements can be found in \cite{Error5}, where Casas and Tr\"{o}ltzsch proved order $O(h)$ in the case of linear-quadratic control problems. Later Casas \cite{linearerror} proved order $o(h)$ for the control problems governed by semilinear elliptic equations and quite general cost functions. In \cite{Error4} R\"{o}sch for the first time proved that the error order is $O(h^{\frac{3}{2}})$ under some special assumptions on the continuous solutions. However, his proof was just done in one dimension. All previous papers were devoted to the full discretization. Recently, a variational discretization concept is introduced by Hinze \cite{Hinze}. More precisely, the state variable and the state equation are discretized, but there is no discretization of the control. He showed that the control error is of order $O(h^2)$. In certain situations, the same convergence order can also be achieved by a special postprocessing procedure, see Meyer and R\"{o}sch \cite{Error6}.

For the study of optimal control problems with sparsity promoting terms, as far as we know, the first paper devoted to this study is published by Stadler \cite{Stadler}, in which structural properties of the control variables were analyzed in the case of the linear-quadratic elliptic optimal control problem. In 2011, a priori and a posteriori error estimates were first given by Wachsmuth and Wachsmuth in \cite{WaWa} for piecewise linear control discretizations, in which the convergence rate is obtained to be of order $O(h)$ under the $L^2$ norm. However, from the point of view of the algorithm, the resulting discretized $L^1$-norm
\begin{eqnarray}\label{equ:discrete norm}
\|u_h\|_{L^1(\Omega_h)}&:=&\int_{\Omega_h}\big{|}{\sum}_{i=1}^{N_h}u_i\phi_i(x)\big{|}\mathrm{d}x,
\end{eqnarray}
does not have a decoupled form with respect to the coefficients $\{u_i\}$, where $\{\phi_i(x)\}$ are the piecewise linear nodal basis functions.
Hence, the authors introduced an alternative discretization of the $L^1$-norm which relies on a nodal quadrature formula
\begin{eqnarray}\label{equ:approximate discrete norm}
\|u_h\|_{L^{1}_h(\Omega)}&:=&{\sum}_{i=1}^{N_h}|u_i|\int_{\Omega_h}\phi_i(x)\mathrm{d}x.
\end{eqnarray}
Obviously, this quadrature incurs an additional error, although the authors proved that this approximation does not change the order of error estimates. In a sequence of papers \cite{CaHerWa1,CaHerWa2}, for the non-convex case governed by a semilinear elliptic equation, Casas et al. proved second-order necessary and sufficient optimality conditions. Using the second-order sufficient optimality conditions, the authors provide error estimates of order $h$ w.r.t. the $L^\infty$ norm for three different choices of the control discretization. It should be pointed out that, for the piecewise linear control discretization case,
a similar approximation technique to the one introduced in \cite{WaWa} is also used for the discretizations of the $L^2$ norm and $L^1$ norm of the control.

Apart from using $L^1$-norm to induce sparsity, Clason and Kunisch in \cite{ClKu1} investigated elliptic control problems with measure-valued controls to promote the sparsity of the control. They discussed the existence and uniqueness of the corresponding dual problems. Subsequently, in 2012, Casas et al in \cite{CaClKu} studied the optimality conditions and provided a priori finite element error estimates for the case of linear-quadratic elliptic control problems with a measure-valued control, in which the control measure was approximated by a linear combination of Dirac measures.

To numerically solve the problem (\ref{eqn:orginal problems}), there are two possible ways.
One is called \emph{First discretize, then optimize}, another approach is called \emph{First optimize, then discretize} \cite{CollisHeink}. Independently of where discretization is located, the resulting finite dimensional equations are quite large. Thus, both of these cases require us to consider proposing an efficient algorithm. In this paper, we focus on the \emph{First discretize, then optimize} approach to solve (\ref{eqn:orginal problems}) and employ the piecewise linear finite elements to discretize (\ref{eqn:orginal problems}).

Next, let us mention some existing numerical methods for solving problem (\ref{eqn:orginal problems}). Since problem (\ref{eqn:orginal problems}) is nonsmooth, thus applying semismooth Newton (SSN) methods is used to be a priority in consideration of their locally superlinear convergence. A special semismooth Newton method with the active set strategy, called the primal-dual active set (PDAS) method is introduced in \cite{BeItKu} for control constrained elliptic optimal control problems. It is proved to have the locally superlinear convergence (see \cite{Ulbrich1,Ulbrich2,HiPiUl} for more details). Mesh-independence results for the SSN method were established in \cite{meshindependent}. Additionally, the authors in \cite{Preconditioning for L1 control} showed that a saddle point system with $2\times 2$ block structure should be solved by employing some Krylov subspace methods with a good preconditioner at each iteration step of the SSN method. However, the $2\times 2$ block linear system is obtained by reducing a $3\times 3$ block linear system with bringing additional computation for linear system involving the mass matrix. Furthermore, the coefficient matrix of the Newton equation would change with every iteration due to the change of the active set. In this case, it is clear that forming a uniform preconditioner, which used to precondition the Krylov subspace methods for solving the Newton equations, is difficult. For a survey of how to precondition saddle point problems, we refer to \cite{Preconditioning for optimal control}.

More importantly, although employing the SSN method can derive the solution with high precision, it is generally known that the total error of utilizing numerical methods to solve PDE constrained problem consists of two parts: discretization error and the iteration error resulted from algorithm of solving the discretized problem. However, the discretization error order for the piecewise linear discretization is $\mathcal{O}(h)$ which accounts for the main part. Thus, algorithms of high precision do not reduce the order of the total error but waste computations. Taking the precision of discretization error into account, employing an efficient first-order algorithms with the aim of solving discretized problems to moderate accuracy is sufficient.

As one may know, for finite dimensional large scale optimization problems, some efficient first-order algorithms, such as iterative soft thresholding algorithms (ISTA) \cite{Blumen}, accelerated proximal gradient (APG)-based method \cite{inexact APG,Beck,Toh}, ADMM \cite{Fazel,SunToh1,SunToh2,SunToh3}, etc, have become the state of the art algorithms. Motivated by the success of these finite dimensional optimization algorithms, Song et al.\cite{iwADMM} proposed an inexact heterogeneous ADMM (ihADMM) for problem (\ref{eqn:orginal problems}). Different from the classical ADMM, the ihADMM adopts two different weighted norms for the augmented term in two subproblems, respectively. Furthermore, the authors also gave theoretical results on the global convergence as well as the iteration complexity results $O(1/k)$.
Recently, thanks to the iteration complexity $O(1/k^2)$, an APG method in function space was proposed to solve (\ref{eqn:orginal problems}) in \cite{FIP}. As we know, the efficiency of the APG method depends on how close the step-length is to the Lipschitz constant. However, in general, choosing an appropriate step-length is difficult since the Lipschitz constant is usually not available analytically. Thus, this disadvantage largely limits the efficiency of APG method.


As far as we know, most of the aforementioned papers are devoted to solving the primal problem. However, when the primal problem (\ref{eqn:orginal problems}) is discretized by the piecewise linear finite elements and directly solved by some algorithms, e.g., SSN, PDAS, ihADMM and APG, as we mentioned above, the resulting discretized $L^1$-norm does not have a decoupled form. Thus the same technique as in (\ref{equ:approximate discrete norm}) should be used, which however will inevitably cause additional error. In this paper, in order to avoid the additional error, we will consider using the duality-based approach for (\ref{eqn:orginal problems}). 
The dual of problem (\ref{eqn:orginal problems}) can be written, in its equivalent minimization form, as
\begin{equation}\label{eqn:dual problem}
\begin{aligned}
\min\ \Phi(\lambda,\mu,p):=&{\frac{1}{2}\|A^*p-y_d\|_{L^2(\Omega)}^2}+ \frac{1}{2\alpha}\|-p+\lambda+\mu\|_{L^2(\Omega)}^2+\langle p, y_r\rangle_{L^2(\Omega)}\\
&+\delta_{\beta B_{\infty}(0)}(\lambda)+\delta^*_{U_{ad}}(\mu)-\frac{1}{2}\|y_d\|_{L^2(\Omega)}^2,\end{aligned}\tag{$\mathrm{D}$}
\end{equation}
where $p\in H^1_0(\Omega)$, $\lambda,\mu\in L^2(\Omega)$, $B_{\infty}(0):=\{\lambda\in L^2(\Omega): \|\lambda\|_{L^\infty(\Omega)}\leq 1\}$, and for any given nonempty, closed convex subset $C$ of $L^2(\Omega)$, $\delta_{C}(\cdot)$ is the indicator function of $C$. Based on the $L^2$-inner product, we define the conjugate of $\delta_{C}(\cdot)$ as follows
\begin{equation*}
 \delta^*_{C}(w^*)=\sup\limits_{w\in C}^{}{\langle w^*,w\rangle}_{L^2(\Omega)}.
\end{equation*}
Although the duality-based approach has been introduced in \cite{ClKu1} for elliptic control problems without control constraints in non-reflexive Banach spaces, the authors did not take advantage of the structure of the dual problem and still used semismooth Newton methods to solve the Moreau-Yosida regularization of the dual problem. In the paper, in terms of the structure of problem (\ref{eqn:dual problem}), we aim to design an algorithm which could efficiently and fast solve the dual problem (\ref{eqn:dual problem}).

By setting $x=(\mu,\lambda,p)$, $x_0=\mu$ and $x_1=\lambda$, it is quite clear that our dual problem (\ref{eqn:dual problem}) belongs to a general class of multi-block convex optimization problems of the form
\begin{equation}\label{equ:general convex problem}
  \min {F(x_0,x):=\varphi(x_0)+\psi(x_1)+\phi(x_0,x)},
\end{equation}
where $x_0\in X_0$, $x=(x_1,...,x_s)\in X:=X_1\times...\times X_s$ and each $X_i$ is a finite dimensional real Euclidean space. The functions $\varphi,\psi$ and $\phi$ are three closed proper convex functions.
Thanks to the structure of (\ref{equ:general convex problem}), in 2015, Chambolle and Pock \cite{Chambolle} proposed the accelerated alternative descent (AAD) algorithm to solve it for this situation that the joint objective function $\phi$ was quadratic . But the disadvantage is that the AAD method does not take the inexactness of the solutions of the associated subproblems into account. As we know, in some case, it is either impossible or extremely expensive to exactly compute the solution of each subproblem even if it is doable, especially at the early stage of the whole process. For example, if a subproblem is equivalent to solving a large-scale or ill-condition linear system, it is a natural idea to use the iterative methods such as some Krylov-based methods. Hence, it is not suitable for the practical application. Subsequenctly, when $\phi$ is a general closed proper convex function and $\arg\min_{x_0}{\varphi(x_0) + \phi(x_0,x)}$ could be computed
exactly, Sun, Toh and Yang \cite{inexact ABCD} proposed an inexact accelerated block coordinate descent (iABCD) method to solve least squares semidefinite programming (LSSDP) via its dual. The basic idea of the iABCD method is firstly applying the Danskin-type theorem to reduce the two block nonsmooth terms into only one block and then using APG method to solve the reduced problem. More importantly, the powerful inexact symmetric Gauss-Seidel (sGS) decomposition technique developed in \cite{SunToh3} is the key for designing the iABCD method. Additionally, the authors proved that the iABCD method has the $O(1/k^2)$ iteration complexity when the subproblems are solved approximately subject to certain inexactness criteria.

However, for the situation the subproblem with respect to block $x_0$ could not be solved exactly,
one could not no longer use Danskin-type theorem to achieve the goal of reducing it into one block nonsmooth term.
To overcome the above bottlenecks, in her PhD thesis \cite[Chapter 3]{CuiYing}, Cui proposed an inexact majorized accelerated block coordinate descent (imABCD) method for solving the following unconstrained convex optimization problems with coupled objective functions
\begin{equation}\label{eqn:model problem1}
\begin{aligned}
\min_{v, w} f(v)+ g(w)+ \phi(v, w).
\end{aligned}
\end{equation}
Under suitable assumptions and certain inexactness criteria, the author can prove that the inexact mABCD method   also enjoys the impressive $O(1/k^2)$ iteration complexity.

In this paper, which is inspired by the success of the iABCD and imABCD methods, we combine their virtues and propose an inexact sGS based majorized ABCD method (called sGS-imABCD) to solve problem (\ref{eqn:dual problem}). The design of this method combines an inexact 2-block majorized ABCD and the recent advances in the inexact sGS technique. Owing to the convergence results of imABCD method which are given in \cite[Chapter 3]{CuiYing}, our proposed algorithm could be proven having the $O(1/k^2)$ iteration complexity as well.

Moreover, some truly implementable inexactness criteria controlling the accuracy of the generated imABCD subproblems are analyzed. Specifically, as shown in Section \ref{sec:5}, because of two nonsmooth subproblems having the closed form solutions, it is easy to see that the main computation of our sGS-imABCD algorithm is in solving $p$-subproblems, which equivalent to solving the $2\times 2$ block saddle point linear system twice at each iteration. It should be pointed out that the coefficient matrix of the saddle point linear system is fixed. To efficiently solve the linear system, a preconditioned GMRES method is used which leads to the rapid convergence and the robustness with respect to the mesh size $h$. More importantly, at first glance, it appears that we would need to solve the linear system twice. In practice, in order to avoid this situation and improve the efficiency of our sGS-imABCD algorithm, we design a strategy to approximate the solution for the second linear system. Thus, when a residual error condition is satisfied, the linear system need only to be solved once instead of twice. We should emphasize that such a saving can be significant, especially in the middle and later stages of the whole algorithm. Thus, in terms of the amount of calculation and the discretized error, our sGS-imABCD algorithm is superior to the semi-smooth Newton method.

As far as we know, we are the first to utilize the duality-based approach and introduce the sGS-imABCD method to solve (\ref{eqn:orginal problems}). In other words, we directly use the sGS-imABCD method to solve problem (\ref{eqn:discretized matrix-vector dual problem}), e.g., the discretized form of the dual problem (\ref{eqn:dual problem}). As already mentioned, one can also apply the ihADMM and APG methods to solve a kind of approximate discretized form (\ref{equ:approx discretized matrix-vector form}) of (\ref{eqn:orginal problems}), where the quadrature technique (\ref{equ:approximate discrete norm}) is used. For the sake of the numerical comparison, we also use our sGS-imABCD method to solve (\ref{eqn:approx discretized matrix-vector dual problem}), e.g., the dual of (\ref{equ:approx discretized matrix-vector form}). As one can see later from the numerical experiments, directly solving (\ref{eqn:discretized matrix-vector dual problem}) can get better discrete error results than that from solving (\ref{eqn:approx discretized matrix-vector dual problem}) and (\ref{equ:approx discretized matrix-vector form}). More importantly, the numerical results also show our sGS-imABCD method is more efficient than the ihADMM and APG methods.

The remainder of the paper is organized as follows. In Section \ref{sec:2}, the first-order optimality conditions for problem (\ref{eqn:orginal problems}) are derived.
 In Section \ref{sec:3}, the finite element approximation is introduced. In Section \ref{sec:4}, we give a review of the inexact sGS technique developed in \cite{SunToh3}, which lays the foundation for further algorithmic developments. In Section \ref{sec:5}, we give a brief sketch of the imABCD \cite[Chapter 3]{CuiYing} and propose our inexact symmetric Gauss-Seidel based majorized ABCD (sGS-imABCD) method. In Section \ref{sec:6}, by comparison with the ihADMM and APG methods, numerical results are given to show the efficiency of our proposed method and confirm the finite element error estimates. Finally, we conclude our paper in Section \ref{sec:7}.

\section{First-order optimality condition}
\label{sec:2}
In this section, we will derive the first-order optimality conditions. First, we analyze the existence and uniqueness of the global solution to problem (\ref{eqn:orginal problems}). Utilizing the Lax-Milgram lemma, we have the following proposition.
\begin{proposition}{\rm\textbf{\cite [Theorem. B.4] {KiSt}}}\label{equ:weak formulation}
Under Assumption {\rm \ref{equ:assumption:1}}, the bilinear form $a(\cdot,\cdot)$ in {\rm (\ref{eqn:bilinear form})} is bounded and $V$-coercive for $V=H^1_0(\Omega)$ and the associate operator $A$ has a bounded inverse. In particular, for every $u \in L^2(\Omega)$ and $y_r\in L^2(\Omega)$, {\rm (\ref{eqn:state equations})} has a unique weak solution $y\in H^1_0(\Omega)$ given by {\rm (\ref{eqn:weak form})}. Furthermore,
\begin{equation}\label{equ:control estimats state}
  \|y\|_{H^1}\leq C (\|u\|_{L^2(\Omega)}+\|y_r\|_{L^2(\Omega)}),
\end{equation}
for a constant $C$ depending only on $a_{ij}$, $c_0$ and $\Omega$.
\end{proposition}
By Proposition \ref{equ:weak formulation} and the strong convexity of the objective function $J(y,u)$ for (\ref{eqn:orginal problems}), it is easy to establish the existence and uniqueness of the solution to (\ref{eqn:orginal problems}). The optimal solution can be characterized by the following Karush-Kuhn-Tucker (KKT) conditions.

\begin{theorem}[{\rm First-Order Optimality Condition}]\label{First-order optimality condition}
Under Assumption \ref{equ:assumption:1}, the couple function ($y^*$, $u^*$) is the optimal solution of {\rm(\ref{eqn:orginal problems})}, if and only if there exists an adjoint state $p^*\in H_0^1(\Omega)$, such that the following conditions hold in the weak sense
\begin{subequations}\label{eqn:KKT}
\begin{eqnarray}
&&\left\{ \begin{aligned}\label{eqn1:KKT}
        Ly^*=u^*+y_r, \quad &\mathrm{in}\ \Omega,\\
        y^*=0,  \quad &\mathrm{on}\ \partial\Omega,
        \end{aligned} \right. \\
&&\left\{ \begin{aligned}\label{eqn2:KKT}
        Lp^*&=y_d-y^*,&  &\mathrm{in}\ \Omega,\\
        p^*&=0,& &\mathrm{on}\ \partial\Omega,
        \end{aligned} \right.\\
&&u^*={\rm\Pi}_{U_{ad}}\left(\frac{1}{\alpha}{\rm{soft}}\left(p^*,\beta\right)\right),\label{eqn3:KKT}
\end{eqnarray}
\end{subequations}
where
\begin{equation*}
  \begin{aligned}
 {\rm\Pi}_{U_{ad}}(v(x))&:=\max\{a,\min\{v(x),b\}\}, \\
 {\rm soft}(v(x),\beta)&:={\rm{sgn}}(v(x))\circ\max(|v(x)|-\beta,0).
  \end{aligned}
\end{equation*}
\end{theorem}

\begin{remark}\label{relation}
From (\ref{eqn3:KKT}), an obvious fact should be pointed out that $|p|<\beta$ implies $u=0$, which also explains that the $L^1$-norm can induce the sparsity property of $u$. Moreover, since $p\in H^1_0(\Omega)$, (\ref{eqn3:KKT}) implies $u\in H^1(\Omega)$. Figure \ref{fig:the relation u and p} shows the relationship between $u$ and $p$
\begin{figure}[H]
\centering
\includegraphics[width=0.40\textwidth]{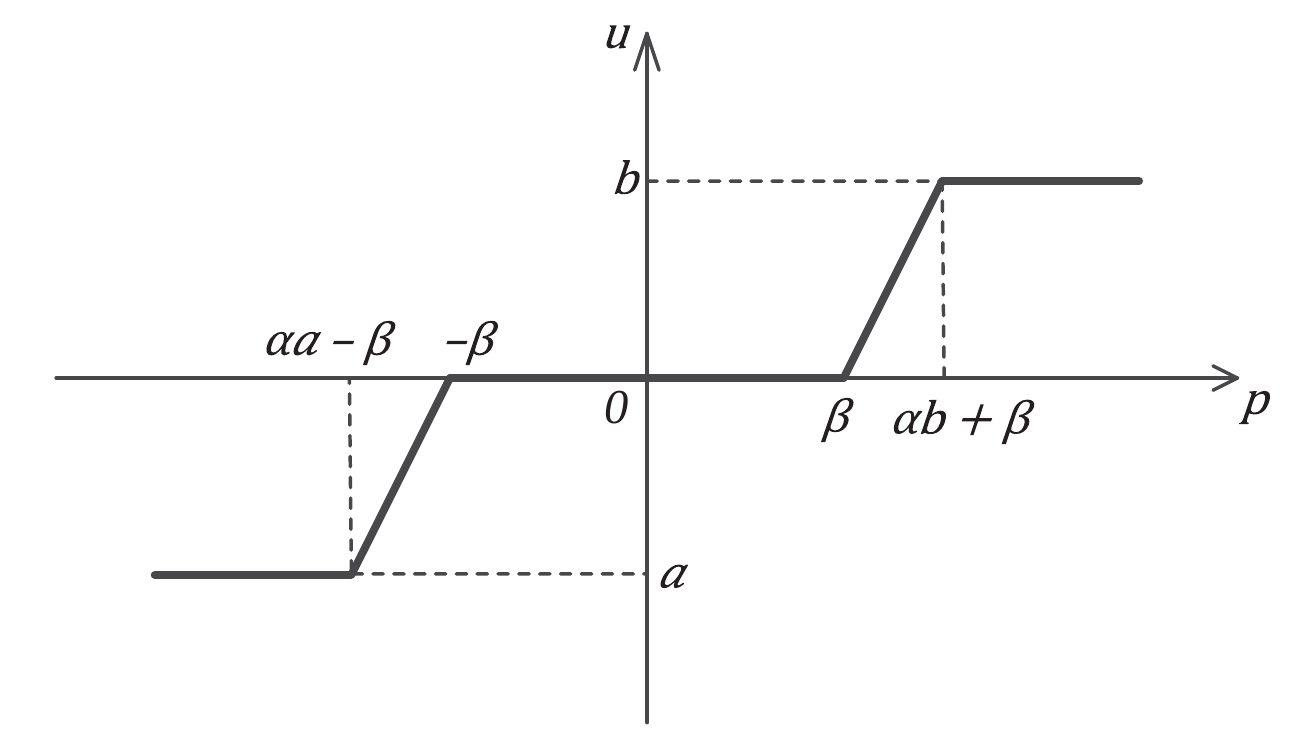}
\caption{ The relationship between $u$ and $p$}\label{fig:the relation u and p}
\end{figure}
\end{remark}
It is obvious that if $\beta$ is sufficiently large, the optimal control would be $u^*_{\beta}=0$. Then we have the following lemma.
\begin{lemma}
  If $\beta\geq\beta_0:=\|(A^{-1})^*(y_d-A^{-1}By_r)\|_{L^\infty{(\Omega)}}$, then the unique solution of problem {\rm(\ref{eqn:orginal problems})} is $(y^*_{\beta},u^*_{\beta})=(A^{-1}By_r,0)$.
\end{lemma}

\section{Finite element approximation}
\label{sec:3}
To numerically solve problem (\ref{eqn:orginal problems}), we consider employing the finite element method, in which the state $y$ and the control $u$ are both discretized by the piecewise linear, globally continuous finite elements.

To this aim, let us fix the assumptions on the discretization by finite elements. We first consider a family of regular and quasi-uniform triangulations $\{\mathcal{T}_h\}_{h>0}$ of $\bar{\Omega}$. For each cell $T\in \mathcal{T}_h$, let us define the diameter of the set $T$ by $\rho_{T}:={\rm diam}\ T$ and define $\sigma_{T}$ to be the diameter of the largest ball contained in $T$. The mesh size of the grid is defined by $h=\max_{T\in \mathcal{T}_h}\rho_{T}$. We suppose that the following regularity assumptions on the triangulation are satisfied which are standard in the context of error estimates.

\begin{assumption}[regular and quasi-uniform triangulations]\label{regular and quasi-uniform triangulations}
There exist two positive constants $\kappa$ and $\tau$ such that
   \begin{equation*}
   \frac{\rho_{T}}{\sigma_{T}}\leq \kappa\ {\rm and}\ \frac{h}{\rho_{T}}\leq \tau
 \end{equation*}
hold for all $T\in \mathcal{T}_h$ and all $h>0$. Moreover, let us define $\bar{\Omega}_h=\bigcup_{T\in \mathcal{T}_h}T$, and let ${\Omega}_h \subset\Omega$ and $\Gamma_h$ denote its interior and its boundary, respectively. In the case that $\Omega$ is a convex polyhedral domain, we have $\Omega=\Omega_h$. In the case $\Omega$ has a $C^{1,1}$- boundary $\Gamma$, we assume that $\bar{\Omega}_h$ is convex and that all boundary vertices of $\bar{\Omega}_h$ are contained in $\Gamma$, such that
\begin{equation*}
  |\Omega\backslash {\Omega}_h|\leq c h^2,
\end{equation*}
where $|\cdot|$ denotes the measure of the set and $c>0$ is a constant.
\end{assumption}

On account of the homogeneous boundary condition of the state equation, we use
\begin{equation*}
  Y_h =\left\{y_h\in C(\bar{\Omega})~\big{|}~y_{h|T}\in \mathcal{P}_1~ {\rm{for\ all}}~ T\in \mathcal{T}_h~ \mathrm{and}~ y_h=0~ \mathrm{in } ~\bar{\Omega}\backslash {\Omega}_h\right\}
\end{equation*}
as the discretized state space, where $\mathcal{P}_1$ denotes the space of polynomials of degree less than or equal to $1$. As mentioned above, we also use the same discretized space to discretize the control $u$, thus we define
\begin{equation*}
   U_h =\left\{u_h\in C(\bar{\Omega})~\big{|}~u_{h|T}\in \mathcal{P}_1~ {\rm{for\ all}}~ T\in \mathcal{T}_h~ \mathrm{and}~ u_h=0~ \mathrm{in } ~\bar{\Omega}\backslash{\Omega}_h\right\}.
\end{equation*}
For a given regular and quasi-uniform triangulation $\mathcal{T}_h$ with nodes $\{x_i\}_{i=1}^{N_h}$, let $\{\phi_i\} _{i=1}^{N_h}$ be a set of nodal basis functions, which span $Y_h$ as well as $U_h$ and satisfy the following properties:
\begin{eqnarray}\label{basic functions properties}
  &&\phi_i \geq 0, \quad
  \|\phi_i\|_{\infty} = 1 \quad \forall i=1,2,...,N_h,
 \quad \sum\limits_{i=1}^{N_h}\phi_i(x)=1\quad {\rm for\ a.a.} \ x\in \Omega_h.
\end{eqnarray}
The elements $u_h\in U_h$ and $y_h\in Y_h$ can be represented in the following forms, respectively,
\begin{equation*}
  u_h=\sum \limits_{i=1}^{N_h}u_i\phi_i,\quad y_h=\sum \limits_{i=1}^{N_h}y_i\phi_i,
\end{equation*}
where $u_h(x_i)=u_i$ and $y_h(x_i)=y_i$. Let $U_{ad,h}$ denote the discretized feasible set, which is defined by
\begin{eqnarray*}
  U_{ad,h}:&=&U_h\cap U_{ad}\\
           &=&\left\{z_h=\sum \limits_{i=1}^{N_h} z_i\phi_i~\big{|}~a\leq z_i\leq b, \forall i=1,...,N_h\right\}\subset U_{ad}.
\end{eqnarray*}
Now, a discretized version of problem (\ref{eqn:orginal problems}) is formulated as follows.
\begin{equation}\label{eqn:discretized problems}
  \left\{ \begin{aligned}
        &\min \limits_{(y_h,u_h)\in Y_h\times U_h}^{}J_h(y_h,u_h)=\frac{1}{2}\|y_h-y_d\|_{L^2(\Omega_h)}^{2}+\frac{\alpha}{2}\|u_h\|_{L^2(\Omega_h)}^{2}+ \beta\|u_h\|_{L^1(\Omega_h)}\\
        &\qquad\quad {\rm{s.t.}}\qquad a(y_h, v_h)=\int_{\Omega}(u_h+y_r)v_h{\rm{d}}x \quad  \forall v_h\in Y_h,  \\
          &\qquad \qquad \qquad~  u_h\in U_{ad,h}.
                          \end{aligned} \right.
 \end{equation}
About the error estimates, we have the following result.
\begin{theorem}{\rm{\textbf{\cite[Proposition 4.3]{WaWa}}}}\label{theorem:error1}
Let us assume that $u^*$ and $u^*_h$ be the optimal control solutions of {\rm(\ref{eqn:orginal problems})} and {\rm(\ref{eqn:discretized problems})}, respectively. Then for every $\alpha_0>0$, $h_0>0$ there exists a constant $C>0$,  such that for all $\alpha\leq\alpha_0$, $h\leq h_0$ the following inequality holds
\begin{equation*}
  \|u-u_h\|_{L^2(\Omega)}\leq C(\alpha^{-1}h+\alpha^{-3/2}h^2),
\end{equation*}
where $C$ is independent of $\alpha, h$.
\end{theorem}

From the perspective of numerical implementation, we introduce the following stiffness and mass matrices
\begin{equation*}
 K_h = \left(a(\phi_i, \phi_j)\right)_{i,j=1}^{N_h},\quad  M_h=\left(\int_{\Omega_h}\phi_i(x)\phi_j(x){\mathrm{d}}x\right)_{i,j=1}^{N_h},
\end{equation*}
and let $y_{r,h}$ and $y_{d,h}$ be the $L^2$-projections of $y_r$ and $y_d$ onto $Y_h$, respectively,
\begin{equation*}
  y_{r,h}=\sum\limits_{i=1}^{N_h}y_r^i\phi_i,\quad y_{d,h}=\sum\limits_{i=1}^{N_h}y_d^i\phi_i.
\end{equation*}
Then, identifying discretized functions with their coefficient vectors, we can rewrite (\ref{eqn:discretized problems}) in the following way:
\begin{equation}\label{equ:discretized matrix-vector form}
\left\{\begin{aligned}
        &\min\limits_{y,u}^{}~~ {J}_h(y,u)= \frac{1}{2}\|y-y_{d}\|_{M_h}^{2}+\frac{\alpha}{2}\|u\|_{M_h}^{2}+ \beta\int_{\Omega_h}|\sum_{i=1}^{N_h}{u_i\phi_i(x)}|~\mathrm{d}x\\
        &\ {\rm{s.t.}}\quad K_hy=M_hu+M_hy_r,\\
        &\ \quad\quad \ a\leq u \leq b.
                          \end{aligned} \right.\tag{${\mathrm{P}}_{h}$}
\end{equation}

It is clear that the discretized $L^1$-norm $\int_{\Omega_h}|\sum_{i=1}^{N_h}{u_i\phi_i(x)}|\mathrm{d}x$ cannot be written as a matrix-vector form and is a coupled form with respect to $\{u_i\}$. Thus, the subgradient $\nu_h\in \partial\|u_h\|_{L^1(\Omega_h)}$ will not belong to a finite-dimensional subspace. Hence, if directly solving (\ref{equ:discretized matrix-vector form}), it is inevitable to bring some difficulties into the numerical calculation. To overcome these difficulties, in \cite{WaWa}, the authors introduced the lumped mass matrix $W_h$
which is a diagonal matrix
\begin{equation*}
  W_h:={\rm{diag}}\left(\int_{\Omega_h}\phi_i(x)\mathrm{d}x\right)_{i=1}^{N_h},
\end{equation*}
and defined an alternative discretization of the $L^1$-norm
\begin{equation}\label{equ:approximal L1}
  \|u_h\|_{L^{1}_h(\Omega)}:=\sum_{i=1}^{N_h}|u_i|\int_{\Omega_h}\phi_i(x)\mathrm{d}x=\|Wu\|_1,
\end{equation}
which is a weighted $l^1$-norm of the coefficients of $u_h$.
More importantly, the following results about the mass matrix $M_h$ and the lumped mass matrix $W_h$ hold.
\begin{proposition}{\rm\textbf{\cite[Table 1]{Wathen}}}\label{eqn:martix properties}
$\forall$ $z\in \mathbb{R}^{N_h}$, the following inequalities hold:
\begin{eqnarray*}
 &&\|z\|^2_{M_h}\leq\|z\|^2_{W_h}\leq c_n\|z\|^2_{M_h}, \quad where \quad c_n=
 \left\{ \begin{aligned}
         &4  \quad if \ n=2 \\
         &5  \quad if \ n=3.
                           \end{aligned} \right.
                           \\
 &&\int_{\Omega_h}|\sum_{i=1}^n{z_i\phi_i(x)}|~\mathrm{d}x\leq\|Wz\|_1.
\end{eqnarray*}
\end{proposition}
Thus, we provide another discretization of problem {\rm(\ref{eqn:orginal problems})}:
\begin{equation}\label{equ:approx discretized matrix-vector form}
\left\{\begin{aligned}
        &\min\limits_{y,u}^{}~~ {J}_h(y,u)= \frac{1}{2}\|y-y_{d}\|_{M_h}^{2}+\frac{\alpha}{2}\|u\|_{M_h}^{2}+ \beta\|W_h u\|_{1}\\
        &\ {\rm{s.t.}}\quad K_hy=M_h (u +y_r),\\
        &\ \quad\quad \ a\leq u \leq b.
                          \end{aligned} \right.\tag{$\widetilde{\mathrm{P}}_{h}$}
\end{equation}
Clearly, the approximation of $L^1$-norm (\ref{equ:approximal L1}) inevitably brings additional error, although it can be proven that this additional error do not disturb the order of error estimates, (see \cite[Corollary 4.6]{WaWa}).

As already mentioned, in this paper, we consider solving problem (\ref{equ:discretized matrix-vector form}) by a duality-based approach. Thus, for the purpose of numerical
implementation, we first give the finite element discretizations of (\ref{eqn:dual problem}) as follows
\begin{equation}\label{eqn:discretized matrix-vector dual problem}
\begin{aligned}
\min\limits_{\mu,\lambda,p\in \mathbb{R}^{N_h}}
\Phi_h(\mu,\lambda,p):=&
\frac{1}{2}\|K_h p-{M_h}y_{d}\|_{M_h^{-1}}^2+ \frac{1}{2\alpha}\|\lambda+ \mu- p\|_{M_h}^2+\langle M_hy_r, p\rangle\\
&+ \delta_{[-\beta,\beta]}(\lambda)+ \delta^*_{[a,b]}({M_h}\mu)-\frac{1}{2}\|y_d\|^2_{M_h}.
\end{aligned}\tag{$\mathrm{D_h}$}
\end{equation}
It is clear that problem (\ref{eqn:discretized matrix-vector dual problem}) is a convex composite minimization problem whose objective is the sum of a coupled quadratic function involving three blocks of variables and two separable non-smooth functions involving only the first and second block, respectively. In the following sections, benefiting from the structure of (\ref{eqn:discretized matrix-vector dual problem}), we aim to propose an efficient and fast algorithm to solve it.

Moreover, for the sake of comparison of numerical experiments, we provide the dual problem of (\ref{equ:approx discretized matrix-vector form}), in its equivalent minimization form, as
\begin{equation}\label{eqn:approx discretized matrix-vector dual problem}
\begin{aligned}
\min\limits_{\mu,\lambda,p\in \mathbb{R}^{N_h}}
\widetilde{\Phi}_h(\mu,\lambda,p):=&
\frac{1}{2}\|K_h p-{M_h}y_{d}\|_{M_h^{-1}}^2+ \frac{1}{2\alpha}\|\lambda+ \mu- p\|_{M_h}^2+\langle M_hy_r, p\rangle\\
&+ \delta_{[-\beta,\beta]}(W_h^{-1}M_h\lambda)+ \delta^*_{[a,b]}({M_h}\mu)-\frac{1}{2}\|y_d\|^2_{M_h}
.\end{aligned}\tag{$\widetilde{\mathrm{D}}_h$}
\end{equation}

\section{An inexact block symmetric Gauss-Seidel iteration}
\label{sec:4}
In this section, we first introduce the symmetric Gauss-Seidel (sGS) technique proposed recently by Li, Sun and Toh \cite{SunToh2}. It is a powerful tool to solve a convex minimization problem whose objective is the sum of a multi-block quadratic function and a non-smooth function involving only the first block, which plays an important role in our subsequent algorithms designs for solving the PDE-constraints optimization problems.

Let $s\geq2$ be a given integer and $\mathcal{X}:=\mathcal{X}_1\times\mathcal{X}_2\times...\times\mathcal{X}_s $ where each $\mathcal{X}_i$ is a real finite dimensional Euclidean space. The sGS technique aims to solve the following unconstrained nonsmooth convex optimization problem approximately
\begin{equation}\label{equ:sGSoptimization}
  \min \phi(x_1)+\frac{1}{2}\langle x, \mathcal{H}x\rangle-\langle r,x\rangle,
\end{equation}
where $x\equiv(x_1,...,x_s)\in \mathcal{X}$ with $x_i\in\mathcal{X}_i$, $i=1,...,s$, $\phi:\mathcal{X}_1\rightarrow (-\infty,+\infty]$ is a closed proper convex function, $\mathcal{H}: \mathcal{X}\rightarrow \mathcal{X}$ is a given self-adjoint positive semidefinite linear operator and $r\equiv(r_1,...,r_s)\in \mathcal{X}$ is a given vector.

For notational convenience, we denote the quadratic function in (\ref{equ:sGSoptimization}) as
\begin{equation}\label{equ:quadratic function}
  h(x):=\frac{1}{2}\langle x, \mathcal{H}x\rangle-\langle r,x\rangle,
\end{equation}
and the block decomposition of the operator $\mathcal{H}$ as
\begin{equation}\label{Hdecomposition}
   \mathcal{H}x:=\left(
    \begin{array}{cccc}
      \mathcal{H}_{11} & \mathcal{H}_{12} & \cdots & \mathcal{H}_{1s} \\
      \mathcal{H}^*_{12} & \mathcal{H}_{22} & \cdots &\mathcal{H}_{2s}  \\
      \vdots & \vdots & \ddots & \vdots \\
      \mathcal{H}^*_{1s} & \mathcal{H}^*_{2s} & \cdots & \mathcal{H}_{ss} \\
    \end{array}
  \right)\left(
                 \begin{array}{c}
                   x_1 \\
                   x_2 \\
                   \vdots \\
                   x_s \\
                 \end{array}
               \right),
\end{equation}
where $\mathcal{H}_{ii}: \mathcal{X}_i\rightarrow \mathcal{X}_i, i=1,...,s$ are self-adjoint positive semidefinite linear operators, $\mathcal{H}_{ij}: \mathcal{X}_j\rightarrow \mathcal{X}_i, i=1,...,s-1, j>i$ are linear maps whose adjoints are given by $\mathcal{H}^*_{ij}$. Here, we assume that $\mathcal{H}_{ii}\succ0, \forall i=1,...,s$. Then, we consider a splitting of $\mathcal{H}$
\begin{equation}\label{equ:splitting}
  \mathcal{H}=\mathcal{D}+\mathcal{U}+\mathcal{U}^*,
\end{equation}
where
\begin{equation}\label{upper triangular}
\mathcal{U}:=\left(
    \begin{array}{cccc}
      0 & \mathcal{H}_{12} & \cdots & \mathcal{H}_{1s} \\
       & \ddots & \cdots &\mathcal{H}_{2s}  \\
      & & \ddots & \mathcal{H}_{(s-1)s} \\
       &  &  & 0 \\
    \end{array}
  \right),
\end{equation}
denotes the strict upper triangular part of $\mathcal{H}$ and $\mathcal{D}:=\mathrm{Diag}(\mathcal{H}_{11},...,\mathcal{H}_{ss})\succ0$ is the diagonal of $\mathcal{H}$. For later discussions, we also define the following self-adjoint positive semidefinite linear operator
\begin{equation}\label{SGSoperator}
  {\rm sGS}(\mathcal{H}):=\mathcal{T}=\mathcal{U}\mathcal{D}^{-1}\mathcal{U}^*.
\end{equation}

For any $x\in \mathcal{X}$, we define
\begin{equation*}
  x_{\leq i}:=(x_1,x_2,...,x_i), \quad x_{\geq i}:=(x_i,x_{i+1},...,x_s),\quad i=0,1,...,s+1,
\end{equation*}
with the convention $x_{\leq 0}=x_{\geq 0}=\emptyset$. Moreover, in order to solve problem (\ref{equ:sGSoptimization}) inexactly, we introduce the following two error tolerance vectors:
\begin{equation*}
  \delta':\equiv(\delta'_1,...,\delta'_s),\quad \delta:\equiv(\delta_1,...,\delta_s),
\end{equation*}
with $\delta'_1=\delta_1$. Define
\begin{equation}\label{sGSerrorterm}
  \Delta(\delta',\delta)=\delta+\mathcal{U}\mathcal{D}^{-1}(\delta-\delta').
\end{equation}
Given $\bar x\in \mathcal{X}$, we consider solving the following problem
\begin{equation}\label{sGSiteration}
 x^{+}:={\arg\min}_x\left\{\phi(x_1)+h(x)+\frac{1}{2}\|x-\bar x\|^2_{\mathcal{T}}-\langle x,\Delta(\delta',\delta)\rangle\right\},
\end{equation}
where $\Delta(\delta',\delta)$ could be regarded as the error term. Then, the following sGS decomposition theorem, which is established by Li, Sun and Toh in \cite{SunToh3}, shows that computing $x^{+}$ in (\ref{sGSiteration}) is equivalent to computing in an inexact block symmetric Gauss-Seidel type sequential updating of the variables $x_1,...,x_s$.
\begin{theorem}\label{sGSTheorem}{\rm\textbf{\cite[Theorem 2.1]{SunToh3}}}
Assume that the self-adjoint linear operators $\mathcal{Q}_{ii}$ are positive definite for all $i=1,...,s$. Then, it holds that
\begin{equation}\label{sGSaddition}
  \mathcal{H}+\mathcal{T}=(\mathcal{D}+\mathcal{U})\mathcal{D}^{-1}(\mathcal{D}+\mathcal{U}^*)\succ0.
\end{equation}
Furthermore, given $\bar x\in \mathcal{X}$, for $i=s,...,2$, suppose we have computed $x'_i\in \mathcal{X}_i$ defined as follows
\begin{equation}\label{backward GS sweep}
  \begin{aligned}
    x'_i:&=\arg\min_{x_i\in\mathcal{X}_i}\phi(\bar x_1)+h(\bar x_{\leq i-1},x_i,x'_{\geq i+1})-\langle \delta'_i,x_i\rangle\\
    &= \mathcal{H}_{ii}^{-1}
    \left(r_i+\delta'_i-\sum\limits^{i-1}_{j=1}\mathcal{H}^*_{ji}
    \bar x_j-\sum\limits^{s}_{j=i+1}\mathcal{H}_{ij}x'_{j}\right),
  \end{aligned}
\end{equation}
then the optimal solution $x^+$ defined by {\rm(\ref{sGSiteration})} can be obtained exactly via
\begin{equation}\label{farward GS sweep}
  \left\{\begin{aligned}
    x^+_1 &=\arg\min_{x_1\in\mathcal{X}_1}\phi(x_1)+h(x_1,x'_{\geq 2})-\langle \delta_1,x_1\rangle,\\
    x^+_i &=\arg\min_{x_i\in\mathcal{X}_i}\phi(x^+_1)+h(x^+_{\leq i-1},x_i,x'_{\geq i+1})-\langle \delta_i,x_i\rangle\\
    &= \mathcal{H}_{ii}^{-1}
    \left(r_i+\delta_i-\sum\limits^{i-1}_{j=1}\mathcal{H}^*_{ji}
    x^+_j-\sum\limits^{s}_{j=i+1}\mathcal{H}_{ij}x'_{j}\right), \quad i=2,...,s.
  \end{aligned}\right.
\end{equation}

\begin{remark}\label{sGSremark}
{\rm(a).} In (\ref{backward GS sweep}) and (\ref{farward GS sweep}), $x'_i$ and $x^+_i$ should be regarded as inexact solutions to the corresponding minimization problems without the linear error terms $\langle \delta'_i,x_i\rangle$ and $\langle \delta_i,x_i\rangle$. Once these approximate solutions
have been computed, they would generate the error vectors $\delta'_i$ and $\delta_i$ as follows:
\begin{equation*}
  \begin{aligned}
  \delta'_i& = \mathcal{H}_{ii}x'_i-
    \left(r_i-\sum\limits^{i-1}_{j=1}\mathcal{H}^*_{ji}\bar x_j-\sum\limits^{s}_{j=i+1}\mathcal{H}_{ij}x'_{j}\right), \quad i=s,...,2,\\
  \delta_1 & \in \partial \phi(x^+_1)+\mathcal{H}_{11}x^+_1-\left(r_1-\sum\limits^{s}_{j=2}\mathcal{H}_{1j}
    x'_j\right),\\
   \delta_i& = \mathcal{H}_{ii}x^+_i-
    \left(r_i-\sum\limits^{i-1}_{j=1}\mathcal{H}^*_{ji}
    x^+_j-\sum\limits^{s}_{j=i+1}\mathcal{H}_{ij}x'_{j}\right), \quad i=2,...,s.
  \end{aligned}
\end{equation*}
With the above known error vectors, we have that $x'_i$ and $x^+_{i}$ are the exact solutions to the minimization problems in {\rm{(\ref{backward GS sweep})}} and {\rm{(\ref{farward GS sweep})}}, respectively.

{\rm(b).} In actual implementations, assuming that for $i=s,...,2$, we have computed $x'_i$ in the backward GS sweep for solving {\rm{(\ref{backward GS sweep})}}, then when solving the subproblems in the forward GS sweep in {\rm{(\ref{farward GS sweep})}} for $i=2,...,s$, we may try to estimate $x^+_i$ by using $x'_i$, and in this case the corresponding error vector $\delta_i$ would be given by
\begin{equation*}
  \delta_i=\delta'_i+\sum\limits^{i-1}_{j=1}\mathcal{H}^*_{ji}(x'_{j}-\bar x_{j}).
\end{equation*}
In practice, we may accept such an approximate solution $x^+_i=x'_i$ for $i=2,...,s$, if the corresponding error vector satisfies an admissible condition such as $\|\delta_i\|\leq c\|\delta'_i\|$ for
some constant $c>1$, say $c = 10$.
\end{remark}
\end{theorem}

In order to estimate the error term $\Delta(\delta',\delta)$ in (\ref{sGSerrorterm}), we have following proposition. 
\begin{proposition}\label{sGSerror}{\rm\textbf{\cite[Proposition 2.1]{SunToh3}}}
Suppose that $\widehat{\mathcal{H}}=\mathcal{H}+\mathcal{T}$ is positive definite. Let $\xi=\|\widehat{\mathcal{H}}^{-1/2}\Delta(\delta',\delta)\|$. It holds that
\begin{equation}\label{equ:sGSerror}
  \xi\leq\|\mathcal{D}^{-1/2}(\delta-\delta')\|+\|\widehat{\mathcal{H}}^{-1/2}\delta'\|.
\end{equation}
\end{proposition}

\section{An inexact majorized accelerated block coordinate descent method for (\ref{eqn:discretized matrix-vector dual problem})}
\label{sec:5}
Obviously, by choosing $v=\mu$ and $w=(\lambda, p)$ and taking
\begin{eqnarray}
  f(v) &=&\delta^*_{[a,b]}({M_h}\mu)\label{f function for Dh},\\
  g(w) &=& \delta_{[-\beta,\beta]}(\lambda)\label{g function for Dh}, \\
  \phi(v, w) &=&\frac{1}{2}\|K_h p-{M_h}y_{d}\|_{M_h^{-1}}^2+ \frac{1}{2\alpha}\|\lambda+ \mu- p\|_{M_h}^2+\langle M_h y_r, p\rangle-\frac{1}{2}\|y_d\|^2_{M_h}\label{phi function for Dh},
\end{eqnarray}
(\ref{eqn:discretized matrix-vector dual problem}) belongs to a general class of unconstrained, multi-block convex optimization problems with coupled objective function, that is
\begin{equation}\label{eqn:model problem}
\begin{aligned}
\min_{v, w} \theta(v,w):= f(v)+ g(w)+ \phi(v, w),
\end{aligned}
\end{equation}
where $f: \mathcal{V}\rightarrow (-\infty, +\infty ]$ and $g: \mathcal{W}\rightarrow  (-\infty, +\infty ]$ are two convex functions (possibly nonsmooth), $\phi: \mathcal{V}\times \mathcal{W}\rightarrow  (-\infty, +\infty ]$ is a smooth convex function, and $\mathcal{V}$, $\mathcal{W}$ are real finite dimensional Hilbert spaces.

\subsection{\bf An imABCD algorithm for general problems (\ref{eqn:model problem})}
It is well known that taking the inexactness of the solutions of associated subproblems into account is important for the numerical implementation. Thus, let us give a brief sketch of the inexact majorized accelerate block coordinate descent (imABCD) method which is proposed by Cui in \cite[Chapter 3]{CuiYing} for the case $\phi$ being a general smooth function. To deal with the general model (\ref{eqn:model problem}), we need some more conditions and assumptions on $\phi$.
\begin{assumption}\label{assumption on differentiable with Lipschitz continuous gradients}
The convex function $\phi: \mathcal{V}\times \mathcal{W}\rightarrow  (-\infty, +\infty ]$ is continuously differentiable with Lipschitz continuous gradient.
\end{assumption}

Let us denote $z:=(v, w)\in \mathcal{V}\times\mathcal{W}$. In \cite[Theorem 2.3]{Generalized Hessian}, Hiriart-Urruty and Nguyen provide a second order Mean-Value Theorem for $\phi$, which states that for any $z'$ and $z$ in $\mathcal{V}\times \mathcal{W}$, there exists $ z''\in [z',z]$ and a self-adjoint positive semidefinite operator $\mathcal{G}\in \partial^2\phi(z'')$ such that
\begin{equation*}
\phi(z)= \phi(z')+ \langle\nabla\phi(z'), z- z'\rangle+ \frac{1}{2}\|z'- z\|_{\mathcal{G}}^2,
\end{equation*}
where $\partial^2\phi(z'')$ denotes the Clarke's generalized Hessian at given $z''$ and $[z',z]$ denotes the
the line segment connecting $z'$ and $z$. Under Assumption \ref{assumption on differentiable with Lipschitz continuous gradients}, it is obvious that there exist two self-adjoint positive semidefinite linear operators $\mathcal{Q}$ and $\widehat{\mathcal{Q}}: \mathcal{V}\times\mathcal{W}\rightarrow \mathcal{V}\times\mathcal{W}$ such that for any
$z\in \mathcal{V}\times\mathcal{W}$,
\begin{equation*}
  \mathcal{Q}\preceq\mathcal{G}\preceq\widehat{\mathcal{Q}}, \quad \forall \ \mathcal{G}\in \partial^2\phi(z).
\end{equation*}
Thus, for any $z, z'\in \mathcal{V}\times\mathcal{W}$, it holds
\begin{equation*}
\phi(z)\geq \phi(z')+ \langle\nabla\phi(z'), z- z'\rangle+ \frac{1}{2}\|z'- z\|_{\mathcal{Q}}^2,
\end{equation*}
and
\begin{equation*}
\phi(z)\leq \hat{\phi}(z; z'):= \phi(z ')+ \langle\nabla\phi(z '), z- z '\rangle+ \frac{1}{2}\|z'- z\|_{\widehat{\mathcal{Q}}}^2.
\end{equation*}
Furthermore, we decompose the operators $\mathcal{Q}$ and $\widehat{\mathcal{Q}}$ into the following block structures
\begin{equation*}
\mathcal{Q}z:=\left(
\begin{array}{cc}
\mathcal{Q}_{11} & \mathcal{Q}_{12}\\
\mathcal{Q}_{12}^* &\mathcal{Q}_{22}
\end{array}
\right)
\left(
\begin{array}{c}
v\\
w
\end{array}
\right),\quad
\widehat{\mathcal{Q}}z:=\left(
\begin{array}{cc}
\widehat{\mathcal{Q}}_{11} & \widehat{\mathcal{Q}}_{12}\\
\widehat{\mathcal{Q}}_{12}^* & \widehat{\mathcal{Q}}_{22}
\end{array}
\right)
\left(
\begin{array}{c}
v\\
w
\end{array}
\right),
\quad\forall z=(v, w)\in \mathcal{U}\times \mathcal{V},
\end{equation*}
and assume $\mathcal{Q}$ and $\widehat{\mathcal{Q}}$ satisfy the following conditions.
\begin{assumption}{\rm\textbf{\cite[Assumption 3.1]{CuiYing}}}\label{assumption majorized}
There exist two self-adjoint positive semidefinite linear operators $\mathcal{D}_1: \mathcal{U}\rightarrow \mathcal{U}$ and $\mathcal{D}_2: \mathcal{V}\rightarrow \mathcal{V}$ such that
\begin{equation*}
\widehat{\mathcal{Q}}:=\mathcal{Q}+ {\rm Diag}(\mathcal{D}_1,\mathcal{D}_2).
\end{equation*}
Furthermore, $\widehat{\mathcal{Q}}$ satisfies that $\widehat{\mathcal{Q}}_{11}\succ 0$ and $\widehat{\mathcal{Q}}_{22}\succ 0$.
\end{assumption}

\begin{remark}\label{choice of Hessian}
It is important to note that Assumption \ref{assumption majorized} is a realistic assumption in practice. For example, when $\phi$ is a quadratic function, we could choose $\mathcal{Q}=\mathcal{G}=\nabla^2\phi$. If we have $\mathcal{Q}_{11}\succ0$ and $\mathcal{Q}_{22}\succ0$, then Assumption \ref{assumption majorized} holds automatically. We should point out that $\phi$ is a quadratic function for many problems in the practical application, such as the SDP relaxation of a binary integer nonconvex quadratic (BIQ) programming, the SDP relaxation for computing lower bounds for quadratic assignment problems (QAPs) and so on, one can refer to \cite{inexact ABCD}. Fortunately, it should be noted that the function $\phi$ defined in (\ref{phi function for Dh}) for our problem (\ref{eqn:discretized matrix-vector dual problem}) is quadratic and thus we can choose $\mathcal{Q}=\nabla^2\phi$.
\end{remark}

We can now present the inexact majorized ABCD algorithm for the general problem (\ref{eqn:model problem}) as follow.

\begin{algorithm}[H]\label{algo 1: imabcd}
  \caption{\textbf{(An inexact majorized ABCD algorithm for (\ref{eqn:model problem}))}}
    \KwIn{$(v^1,w^1)=(\tilde{v}^0, \tilde{w}^0)\in \textrm{dom} (f)\times \textrm{dom}(g)$. Let $\{\epsilon_k\}$ be a summable sequence of nonnegative numbers, and set $t_1=1$, $k=1$.}
  \KwOut{$ (\tilde{v}^k, \tilde{w}^k)$}
  Iterate until convergence:
\begin{description}
  \item[\bf Step 1] Choose error tolerance $\delta_v^k\in \mathcal{U}, \delta_w^k\in \mathcal{V}$ such that
  \begin{equation*}
  \max\{\delta_v^k, \delta_w^k\}\leq \epsilon_k.
  \end{equation*}
Compute
  \begin{equation*}
  \left\{
  \begin{aligned}
    &\tilde{v}^{k}
        =\arg\min_{v \in \mathcal{V}}\{f(v)+ \hat{\phi}(v,w^k; v^k,w^k)- \langle\delta_v^k, v \rangle\},\\
    &\tilde{w}^{k}
    	=\arg\min_{w \in \mathcal{W}}\{g(w)+ \hat{\phi}(\tilde{v}^k,w; v^k,w^k)- \langle\delta_w^k, w \rangle\}.\\
  \end{aligned}
  \right.
  \end{equation*}
  \item[\bf Step 2] Set $t_{k+1}=\frac{1+\sqrt{1+4t_k^2}}{2}$ and $\beta_k=\frac{t_k-1}{t_{k+1}}$, compute
  \begin{equation*}
                v^{k+1}=\tilde{v}^{k}+ \beta_{k}(\tilde{v}^{k}-\tilde{v}^{k-1}), \quad w^{k+1}= \tilde{w}^{k}+ \beta_{k}(\tilde{w}^{k}-\tilde{w}^{k-1}).
  \end{equation*}
\end{description}
\end{algorithm}

Here we state the convergence result without proving. For the detailed proof, one could see \cite[Chapter 3]{CuiYing}. This theorem builds a solid foundation for our our subsequent proposed algorithm.
\begin{theorem}{\rm\textbf{\cite[Theorem 3.2]{CuiYing}}}\label{imABCD convergence}
Suppose that Assumption {\rm\ref{assumption majorized}} holds and the solution set $\Omega$ of the problem {\rm(\ref{eqn:model problem})} is non-empty. Let
$z^*=(v^*,w^*)\in \Omega$. Assume that $\sum\limits_{k=1}^{\infty}k\epsilon_k< \infty$. Then the sequence $\{\tilde{z}^k\}:=\{(\tilde{v}^k,\tilde{w}^k)\}$ generated by the Algorithm {\rm\ref{algo 1: imabcd}} satisfies that
\begin{equation*}
\theta(\tilde{z}^k)- \theta(z^*)\leq \frac{2\|\tilde{z}^0- z^*\|_{\mathcal{S}}^2+ c_0}{(k+1)^2}, \quad \forall k\geq 1,
\end{equation*}
where $c_0$ is a constant number and $\mathcal{S}:={\rm{Diag}}(\mathcal{D}_1,\mathcal{D}_2+\mathcal{Q}_{22})$.
\end{theorem}

\subsection{\bf A sGS-imABCD algorithm for (\ref{eqn:discretized matrix-vector dual problem})}
Now, we can apply Algorithm \ref{algo 1: imabcd} to our problem (\ref{eqn:discretized matrix-vector dual problem}), where $\mu$ is taken as one block, and $(\lambda,p)$ are taken as the other one. Let us denote $z=(\mu, \lambda,p)$. Since $\phi$ defined in (\ref{phi function for Dh}) for (\ref{eqn:discretized matrix-vector dual problem}) is quadratic, we can take
\begin{equation}\label{Hessen Matrix for Dh}
\mathcal{Q}:=
\frac{1}{\alpha}\left(
\begin{array}{ccc}
M_h &\quad M_h & -M_h\\
M_h & \quad  M_h& -M_h\\
-M_h & \quad-M_h&\quad M_h+ \alpha K_h M_h^{-1}K_h
\end{array}
\right),
\end{equation}
where
\begin{equation*}
\mathcal{Q}_{11}:= \frac{1}{\alpha} M_h,\quad
\mathcal{Q}_{22}:=
\frac{1}{\alpha}\left(
\begin{array}{cc}
M_h& \quad-M_h\\
-M_h & \quad M_h+\alpha K_h M_h^{-1}K_h
\end{array}
\right).
\end{equation*}
Additionally, we assume that there exists two self-adjoint positive semidefinite operators $\mathcal{D}_1$ and $\mathcal{D}_2$, such that Assumation \ref{assumption majorized} holds. It implies that we should majorize $\phi(\mu, \lambda,p)$ at $z'=(\mu',\lambda', p')$ as
\begin{equation}\label{majorized function}
\begin{aligned}
 \phi(z) \leq \hat{\phi}(z;z')=& \phi(z)+\frac{1}{2}\|\mu-\mu'\|^2_{\mathcal{D}_1}+\frac{1}{2}\left\|\left(\begin{array}{c}
             \lambda \\
             p
             \end{array}\right)
-\left(\begin{array}{c}
\lambda' \\
p'
\end{array}\right)\right\|^2_{\mathcal{D}_2}.
 \end{aligned}
\end{equation}
Thus, the framework of imABCD for (\ref{eqn:discretized matrix-vector dual problem}) is given below:

\begin{algorithm}[H]\label{algo1:imABCD algorithm for (Dh)}
  \caption{\textbf{(imABCD algorithm for (\ref{eqn:discretized matrix-vector dual problem}))}}
 \KwIn{$(\mu^1, {\lambda}^1, {p}^1)=(\tilde{\mu}^0, \tilde{\lambda}^0, \tilde{p}^0)\in {\rm dom} (\delta^*_{[a,b]})\times [-\beta,\beta]\times \mathbb{R}^{N_h}$. Set $k= 1, t_1= 1.$}

 \KwOut{$ (\tilde{\mu}^k, \tilde{\lambda}^k, \tilde{p}^k)$}

 Iterate until convergence
\begin{description}
  \item[\bf Step 1] Compute 
\begin{eqnarray*}
  \tilde{\mu}^{k}&=&\arg\min\delta^*_{[a,b]}(M_h\mu)+\phi(\mu,\lambda^k,p^k)+\frac{1}{2}\|\mu-\mu^k\|^2_{\mathcal{D}_1} -\langle\delta_{\mu}^k,\mu\rangle,
 \\
(\tilde{\lambda}^{k},\tilde{p}^{k})&=&\arg\min\delta_{[-\beta,\beta]}(\lambda)+\phi(\tilde{\mu}^k,\lambda,p)+\frac{1}{2}\left\|\left(\begin{array}{c}
             \lambda \\
             p
             \end{array}\right)
-\left(\begin{array}{c}
\lambda^k \\
p^k
\end{array}\right)\right\|^2_{\mathcal{D}_2}-\langle\delta_{\lambda}^k,\lambda\rangle-\langle\delta_{p}^k,p\rangle.
\end{eqnarray*}
  \item[\bf Step 2] Set $t_{k+1}=\frac{1+\sqrt{1+4t_k^2}}{2}$ and $\beta_k=\frac{t_k-1}{t_{k+1}}$, compute
\begin{eqnarray*}
\mu^{k+1}=\tilde{\mu}^{k}+ \beta_{k}(\tilde{\mu}^{k}-\tilde{\mu}^{k-1}),\quad
 p^{k+1}=\tilde{p}^{k}+ \beta_{k}(\tilde{p}^{k}-\tilde{p}^{k-1}), \quad
 \lambda^{k+1}= \tilde{\lambda}^{k}+ \beta_{k}(\tilde{\lambda}^{k}-\tilde{\lambda}^{k-1}).
\end{eqnarray*}
\end{description}
\end{algorithm}

Next, another key issue should be considered is how to choose the operators $\mathcal{D}_1$ and $\mathcal{D}_2$. As we know, choosing the appropriate and effective operators $\mathcal{D}_1$ and $\mathcal{D}_2$ is an important thing from the perspective of both theory analysis and numerical implementation. Note that for numerical efficiency, the general principle is that both $\mathcal{D}_1$ and $\mathcal{D}_2$ should be chosen as small as possible such that $\tilde{\mu}^{k}$ and $(\tilde{\lambda}^{k},\tilde{p}^{k})$ could take larger step-lengths while the corresponding subproblems still could be solved relatively easily.

First, for the proximal term $\frac{1}{2}\|\mu-\mu^k\|^2_{\mathcal{D}_1}$, in order to make the subproblem of the block $\mu$ having a analytical solution, and from Proposition (\ref{eqn:martix properties}), we choose
\begin{equation*}
  \mathcal{D}_1:=\frac{1}{\alpha}c_n M_hW_h^{-1}M_h-\frac{1}{\alpha}M_h,\quad {\rm where} \  c_n =
 \left\{ \begin{aligned}
         &4  \quad if \ n=2, \\
         &5  \quad if \ n=3.
                           \end{aligned} \right.
\end{equation*}
For more details, one can see Subsection \ref{subsection5.1}.

Next, we will focus on how to choose the operator $\mathcal{D}_2$. If we ignore the proximal term $\frac{1}{2}\left\|\left(\begin{array}{c}
				\lambda \\
				p
			\end{array}\right)
			-\left(\begin{array}{c}
				\lambda^k \\
				p^k
			\end{array}\right)\right\|^2_{\mathcal{D}_2}$ and the error terms, it is obvious that the subproblem of the block $(\lambda,p)$ belongs to the form (\ref{equ:sGSoptimization}), which can be rewritten
as:
\begin{equation}\label{sGS subproblem for Dh}
  \min \delta_{[-\beta,\beta]}(\lambda)+\frac{1}{2}\langle \left(\begin{array}{c}
                                                             \lambda \\
                                                             p
                                                           \end{array}\right)
  ,\mathcal{H}\left(\begin{array}{c}
                                                             \lambda \\
                                                             p
                                                           \end{array}\right)\rangle-\langle r,\left(\begin{array}{c}
                                                             \lambda \\
                                                             p
                                                           \end{array}\right)\rangle,
\end{equation}
where $\mathcal{H}=
\mathcal{Q}_{22}=
\frac{1}{\alpha}\left(
\begin{array}{cc}
M_h& \quad-M_h\\
-M_h & \quad M_h+\alpha K_h M_h^{-1}K_h
\end{array}
\right)$ and $r=\left(\begin{array}{c}
                                                             \frac{1}{\alpha}M_h\tilde{\mu}^k \\
                                                             M_hy_r-K_hy_d-\frac{1}{\alpha}M_h\tilde{\mu}^k
                                                           \end{array}\right)$.
Since the objective function of (\ref{sGS subproblem for Dh}) is the sum of a two-block quadratic function and a non-smooth function involving only the first block, thus the inexact sGS technique, which is introduced in Section \ref{sec:4}, can be used to solve (\ref{sGS subproblem for Dh}) . To achieve our goal, we choose
\begin{equation*}
  \mathcal{\widetilde{D}}_2={\rm sGS}(\mathcal{Q}_{22})=\frac{1}{\alpha}\left(
                                      \begin{array}{cc}
                                        M_h(M_h+\alpha K_h M_h^{-1}K_h)^{-1}M_h &  \quad0\\
                                        0 &  \quad 0\\
                                      \end{array}
                                    \right).
\end{equation*}
Then according to Theorem \ref{sGSTheorem}, we can solve the $(\lambda,p)$-subproblem by the following procedure
\begin{equation}
\left\{\begin{aligned}\label{sGS procedure of lambda p}
 \hat{p}^{k}&=\arg\min\frac{1}{2}\|K_h p-{M_h}y_{d}\|_{M_h^{-1}}^2+ \frac{1}{2\alpha}\|p-\lambda^k-\tilde{\mu}^k+\alpha y_r\|_{M_h}^2- \langle\hat{\delta}^k_p, p\rangle,\\
 \tilde{\lambda}^{k}
        &=\arg\min\frac{1}{2\alpha}\|\lambda-(\hat{p}^{k}-\tilde{\mu}^k)\|_{M_h}^2+\delta_{[-\beta,\beta]}(\lambda),\\
\tilde{p}^{k}&=\arg\min\frac{1}{2}\|K_h p-{M_h}y_{d}\|_{M_h^{-1}}^2+ \frac{1}{2\alpha}\|p-\tilde{\lambda}^k-\tilde{\mu}^k+\alpha y_r\|_{M_h}^2- \langle\delta^k_p, p\rangle.
\end{aligned}\right.
\end{equation}
However, it is easy to see that the $\lambda$-subproblem is coupled about the variable $\lambda$ since the mass matrix $M_h$ is not diagonal, thus there is no a closed form solution for $\lambda$. To overcome this difficulty, we can take advantage of the relationship between the mass matrix $M_h$ and the lumped mass matrix $W_h$ and add a proximal term $\frac{1}{2\alpha}\|\lambda- \lambda^{k}\|_{W_h-M_h}^2$ to the $\lambda$-subproblem. Fortunately, we have

\begin{equation*}
    {\rm sGS}(\mathcal{Q}_{22})={\rm sGS}\left(\mathcal{Q}_{22}+\frac{1}{\alpha}\left[
                                      \begin{array}{cc}
                                        W_h-M_h &\quad0  \\
                                        0 & \quad0 \\
                                      \end{array}
                                    \right]\right),
\end{equation*}
which implies that the proximal term $\frac{1}{2\alpha}\|\lambda- \lambda^{k}\|_{W_h-M_h}^2$ has no influence on the sGS technique. Thus, we can choose $\mathcal{D}_2$ as follows
\begin{equation*}
    \mathcal{D}_{2}={\rm sGS}(\mathcal{Q}_{22})+\frac{1}{\alpha}\left(
                                      \begin{array}{cc}
                                        W_h-M_h &\quad0  \\
                                        0 & \quad0 \\
                                      \end{array}
                                    \right).
\end{equation*}
Based on the choice of $\mathcal{D}_1$ and $\mathcal{D}_2$, 
we get the majorized Hessian matrix $\widehat{\mathcal{Q}}$ as follows
\begin{equation}\label{eq:majorized Hessen matrix}
\widehat{\mathcal{Q}}
=Q+ \frac{1}{\alpha}
\left(
\begin{array}{ccc}
c_n M_hW_h^{-1}M_h-M_h& 0 &\quad0\\
 0&M_h(M_h+\alpha K_h M_h^{-1}K_h)^{-1}M_h+W_h-M_h &\quad0\\
 0& 0&\quad 0
\end{array}
\right).
\end{equation}
Then, according to the choice of $\mathcal{D}_{1}$ and $\mathcal{D}_{2}$,  we give the detailed framework of our inexact sGS based majorized ABCD method (called sGS-imABCD) for (\ref{eqn:discretized matrix-vector dual problem}) as follows.

\begin{algorithm}[H]\label{algo1:Full inexact ABCD algorithm for (Dh)}
  \caption{\textbf{(sGS-imABCD algorithm for (\ref{eqn:discretized matrix-vector dual problem}))}}
 \KwIn{$(\mu^1, {\lambda}^1, {p}^1)=(\tilde{\mu}^0, \tilde{\lambda}^0, \tilde{p}^0)\in {\rm dom} (\delta^*_{[a,b]})\times [-\beta,\beta]\times \mathbb{R}^{N_h}$. Let $\{\epsilon_k\}$ be a nonincreasing sequence of nonnegative numbers such that $\sum\limits_{k=1}^{\infty}k\epsilon_k< \infty$. Set $k= 1, t_1= 1.$}

 \KwOut{$ (\tilde{\mu}^k, \tilde{\lambda}^k, \tilde{p}^k)$}

 Iterate until convergence
\begin{description}
  \item[\bf Step 1] Choose error tolerance $\delta_{\mu}^k, \hat\delta_{p}^k, \delta_p^k$ such that
  \begin{equation*}
\max\{\|{\delta}_{\mu}^k|\|,\|\hat{\delta}_p^k|\|,\|\delta_p^k|\|\}\leq \epsilon_k.
  \end{equation*}
\begin{itemize}
\item[]Compute 
\begin{eqnarray*}
  \tilde{\mu}^{k}&=&\arg\min\frac{1}{2\alpha}\|\mu-(p^k-\lambda^k)\|_{M_h}^2+\delta^*_{[a,b]}(M_h\mu)+\frac{1}{2}\|\mu-\mu^k\|^2_{\mathcal{D}_1}-\langle\delta^k_{\mu}, \mu\rangle,
  \\
 \hat{p}^{k}&=&\arg\min\frac{1}{2}\|K_h p-{M_h}y_{d}\|_{M_h^{-1}}^2+ \frac{1}{2\alpha}\|p-\lambda^k-\tilde{\mu}^k+\alpha y_r\|_{M_h}^2- \langle\hat{\delta}^k_p, p\rangle,\\
 \\
 \tilde{\lambda}^{k}
        &=&\arg\min\frac{1}{2\alpha}\|\lambda-(\hat{p}^{k}-\tilde{\mu}^k)\|_{M_h}^2+\delta_{[-\beta,\beta]}(\lambda)+\frac{1}{2\alpha}\|\lambda- \lambda^{k}\|_{W_h-M_h}^2,\\
 \\
\tilde{p}^{k}&=&\arg\min\frac{1}{2}\|K_h p-{M_h}y_{d}\|_{M_h^{-1}}^2+ \frac{1}{2\alpha}\|p-\tilde{\lambda}^k-\tilde{\mu}^k+\alpha y_r\|_{M_h}^2- \langle\delta^k_p, p\rangle.
\end{eqnarray*}
\end{itemize}
  \item[\bf Step 2] Set $t_{k+1}=\frac{1+\sqrt{1+4t_k^2}}{2}$ and $\beta_k=\frac{t_k-1}{t_{k+1}}$, compute
\begin{eqnarray*}
\mu^{k+1}=\tilde{\mu}^{k}+ \beta_{k}(\tilde{\mu}^{k}-\tilde{\mu}^{k-1}),\quad
 p^{k+1}=\tilde{p}^{k}+ \beta_{k}(\tilde{p}^{k}-\tilde{p}^{k-1}), \quad
 \lambda^{k+1}= \tilde{\lambda}^{k}+ \beta_{k}(\tilde{\lambda}^{k}-\tilde{\lambda}^{k-1}).
\end{eqnarray*}
\end{description}
\end{algorithm}

Based on Theorem \ref{imABCD convergence}, we can show our Algorithm \ref{algo1:Full inexact ABCD algorithm for (Dh)} (sGS-imABCD) also has the following $O(1/k^2)$ iteration complexity.

\begin{theorem}\label{sGS-imABCD convergence}
Assume that $\sum\limits_{i=k}^{\infty}k\epsilon_k< \infty$. Let $\{\tilde{z}^k\}:=\{(\tilde{\mu}^k,\tilde{\lambda}^k,\tilde{p}^k)\}$ be the sequence generated by the Algorithm \ref{algo1:Full inexact ABCD algorithm for (Dh)}. Then we have
\begin{equation*}
\Phi_h(\tilde{z}^k)- \Phi_h(z^*)\leq \frac{2\|\tilde{z}^0- z^*\|_{\mathcal{S}}^2+ c_0}{(k+1)^2}, \; \forall k\geq 1,
\end{equation*}
where $c_0$ is a constant number, $\mathcal{S}:={\rm{Diag}}(\mathcal{D}_1,\mathcal{D}_2+\mathcal{Q}_{22})$, and $\Phi_h(\cdot)$ is the objective function of the dual problem {\rm(\ref{eqn:discretized matrix-vector dual problem})}.
\end{theorem}
\proof
By Proposition \ref{eqn:martix properties}, we know that
$c_n M_hW_h^{-1}M_h-M_h \succ0$, $M_h(M_h+\alpha K_h M_h^{-1}K_h)^{-1}M_h \succ0$, $W_h-M_h \succ 0$. Moreover, since stiffness and mass matrices are symmetric positive definite matrices, it is noticed that Assumption \ref{assumption majorized} is valid for our $\widehat{\mathcal{Q}}$ which is defined in (\ref{eq:majorized Hessen matrix}). Thus, according to Theorem \ref{imABCD convergence}, we can establish the convergence of Algorithm \ref{algo1:Full inexact ABCD algorithm for (Dh)}.
\endproof
\begin{remark}
Let $\tau_h=2\|\tilde{z}^0- z^*\|_{\mathcal{S}}^2+ c_0$. It is obvious that $\tau_h$ is independent of the parameter $\beta$, whereas it depends on the parameter $\alpha$ and will increase with the decrease of $\alpha$.
\end{remark}

\subsection{\textbf {Numerical computation of the block $\mu$ and $\lambda$ subproblems}}\label{subsection5.1}
For the first subproblem of Algorithm {\rm\ref{algo1:Full inexact ABCD algorithm for (Dh)}} in $k$th iteration, at first glance, there is no closed form solution for the variable $\mu$. However, if we carefully check the subproblems with respect to the variables $p$ and $\lambda$, it is easy to see that we only need the value $M_h\mu$ instead of $\mu$. Thus, let us denote $\xi=M_h\mu$, then solving the subproblem about the variable $\mu$ can be translate to solving the following subproblem
\begin{equation}\label{subproblem-z}
  \begin{aligned}
\tilde{\xi}^{k}
     &=\arg\min\frac{1}{2\alpha}\|\xi-M_h(p^k-\lambda^k)\|_{M_h^{-1}}^2+\delta^*_{[a,b]}(\xi)+\frac{1}{2\alpha}\|\xi-\xi^k\|^2_{c_n W_h^{-1}-M_h^{-1}}\\
     &=\arg\min\frac{1}{2\alpha}\|\xi-(\xi^k+\frac{1}{c_n}W_h(M_h(p^k-\lambda^k)-M_h^{-1}\xi^k))\|_{c_n W_h^{-1}}^2+\delta^*_{[a,b]}(\xi).\\
  \end{aligned}
\end{equation}
To solve (\ref{subproblem-z}), we first introduce the proximal mapping $\textrm{prox}^{f}_{\mathcal{M}}(\cdot)$ with respect to a self-adjoint positive definite linear operator $\mathcal{M}$, which is defined as
\begin{equation}\label{proximal mapping}
  \textrm{prox}^{f}_{\mathcal{M}}(x)=\arg\min\{f(z)+\frac{1}{2}\|z-x\|^2_{\mathcal{M}}\}, \quad \forall x\in \mathcal{\mathcal{X}},
\end{equation}
where $f$ is a closed proper convex function $f$ and $\mathcal{X}$ is a finite-dimensional real Euclidean space.

For the proximal mapping, we have the following Moreau identity which is shown in \cite[Proposition 2.4]{twophaseALM}:
\begin{equation}\label{equ:Moreau identity}
  x=\textrm{prox}^{f}_{\mathcal{M}}(x)+\mathcal{M}^{-1}\textrm{prox}^{f^*}_{\mathcal{M}^{-1}}(\mathcal{M}x),
\end{equation}
where $f^*$ is the conjugate function of $f$. Thus, making use of the Moreau identity (\ref{equ:Moreau identity}), we can derive
\begin{equation}\label{closed form-z}
\tilde{\xi}^{k}=\vartheta^k-\frac{\alpha}{c_n}W_h{\rm\Pi}_{[a,b]}(\frac{c_n}{\alpha}W_h^{-1}\vartheta^k).
\end{equation}
where
\begin{eqnarray*}
  \vartheta^k
  &=& \xi^k+\frac{1}{c_n}W_h(M_h(p^k-\lambda^k)-\mu^k).
\end{eqnarray*}
This means that the subproblem about $\mu$ has a closed form solution. And this is also the important reason why we choose the proximal term $\frac{1}{2\alpha}\|\mu-\mu^k\|^2_{c_n M_hW_h^{-1}M_h-M_h}$ for $\mu$.

When computing the primal variables $y$ and $u$, we still require $\mu$. Then we can compute $\tilde{\mu}^k$ by $\tilde{\mu}^k= M_h^{-1}\tilde{\xi}^{k}$. Based on the eigenvalues bounds for the mass matrix given in \cite{Wathen}, we suggest that using a fix number steps of Chebyshev semi-iteration to represent approximation to $M_h^{-1}$ is an appropriate choice. For more details on the Chebyshev semi-iteration method we refer to {\rm\cite{ReDoWa,chebysev semi-iteration}}. In actual numerical implementations, we use 20 steps of Chebyshev semi-iteration and set the error tolerance to be $10^{-12}$, which could guarantee the error vector $\|{\delta}_{\mu}^k\|_{2} \leq {\epsilon_k}$.

For the block $\lambda$, since the lumped mass matrix $W_h$ is a diagonal positive definite matrix, we can easily derive that
\begin{equation*}
  \tilde{\lambda}^{k}={\rm\Pi}_{[-\beta,\beta]}(s^k),
\end{equation*}
where $s^k=\lambda^k+W_h^{-1}M_h(\hat{p}^k-\tilde{\mu}^k-\lambda^k)$.

\subsection{\textbf {An efficient iteration method and preconditioner for the block $\hat{p}$ subproblem}}\label{subsection5.2}
As we know, the main computation of our sGS-imABCD algorithm is in solving $p$-subproblems. Thus, it is crucial to improve the efficiency of ous sGS-imABCD algorithm in employing an fast strategy to solve $p$-subproblems. For the $\hat{p}^{k}$-subproblem, if we ignore the error vector $\hat{\delta}_{p}^k$, it is obvious to see that solving the subproblem is equivalent to solving the following system:
\begin{equation}\label{discrete optimal condition}
K_hM_h^{-1}(K_h\hat{p}^k-M_hy_d)+\frac{1}{\alpha}M_h(\hat{p}^k-\lambda^k-\tilde{\mu}^k+\alpha y_r)=0.
\end{equation}
Since $K_h{p}=M_h(y_d-y)$, then (\ref{discrete optimal condition}) can be rewritten as:
\begin{equation}\label{equ:discretized saddle point problems2}
\mathcal{A}w^{k+1}\equiv\left[
  \begin{array}{cc}
    \frac{1}{\alpha}M_h & \quad-K_h \\
    K_h & \quad M_h\\
  \end{array}
\right]\left[
         \begin{array}{c}
           \hat{p}^{k} \\
           \hat{y}^{k} \\
         \end{array}
       \right]=\left[
                 \begin{array}{c}
                  \frac{1}{\alpha}M_h(\lambda^k+\tilde{\mu}^k-\alpha y_r) \\
                   M_hy_d\\
                 \end{array}
               \right].
\end{equation}
Clearly, the linear system {\rm(\ref{equ:discretized saddle point problems2})} is a special case of the generalized saddle-point system, thus some Krylov-based methods could be employed to inexactly solve the linear system by constructing a good preconditioner. Here, the preconditioned variant of modified hermitian and skew-hermitian splitting {\rm(PMHSS)} preconditioner
\begin{equation*}
  \mathcal{P_{HSS}}=\frac{1}{\alpha}\left[
    \begin{array}{lcc}
      I & -\sqrt{\alpha}I \\
      \sqrt{\alpha}I & \alpha I \\
    \end{array}
  \right]\left[
           \begin{array}{lcc}
             M_h+\sqrt{\alpha}K_h & 0 \\
             0 & M_h+\sqrt{\alpha}K_h \\
           \end{array}
         \right],
\end{equation*}
which is introduced in {\rm \cite{Bai}}, is employed to precondition the generalized minimal residual {\rm (GMRES)} method to solve {\rm(\ref{equ:discretized saddle point problems2})}.
About the spectral properties of the preconditioned matrix $\mathcal{P_{HSS}}^{-1}\mathcal{A}$, we introduce the following theorem, see {\rm \cite[Theorem 2.3]{Bai}} for more details.
\begin{theorem}\label{spectral properties}
When $\mathcal{P_{HSS}}$ is used to precondition the matrix $\mathcal{A}$, the eigenvalues of the preconditioned matrix $\mathcal{P_{HSS}}^{-1}\mathcal{A}$ are contained within the complex disk centred at $1$ with radius $\frac{\sqrt{2}}{2}$. Moreover, the matrix $\mathcal{P_{HSS}}^{-1}\mathcal{A}$ is diagonalizable.
\end{theorem}
It would be crucial to pointed out that the reason we prefer the PMHSS-preconditioned GMRES method is because it shows $h$- and $\alpha$-independent convergence properties, see the numerical results in {\rm \cite{Bai}} for more details.

In actual implementations, the action of the preconditioning matrix, when used to
precondition the \emph{Krylov} subspace methods, is realized through solving a sequence of generalized residual equations of the form
 \begin{equation*}
   \mathcal{P_{HSS}}v=r,
 \end{equation*}
where $r=(r_a; r_b)\in \mathbb{R}^{2N_h}$, with $r_a, r_b\in \mathbb{R}^{N_h}$, represents the current residual vector, while $v=(v_a; v_b)\in \mathbb{R}^{2N_h}$, with $v_a, v_b\in \mathbb{R}^{N_h}$, represents the generalized residual vector. By making use of the structure of the matrix $\mathcal{P_{HSS}}$, we obtain the following procedure for computing the vector $v$
\begin{algorithm}\label{Algorithm for Precondition}
  \caption{\textbf {Numerical implementation of $\mathcal{P_{HSS}}$ }}
\begin{itemize}
  \item[] {\bf Step 1}. compute $\hat{r}_a$ and $\hat{r}_b$
  \begin{equation*}
    \hat{r}_a=1/2(\alpha r_a+\sqrt{\alpha}r_b), \quad\hat{r}_b=1/2(r_b-\sqrt{\alpha}r_a).
  \end{equation*}
  \item[] {\bf Step 2}. compute $v_a$ and $v_b$ by solving the following linear systerms
  \begin{equation*}
    (M_h+\sqrt{\alpha}K_h)v_a=\hat {r}_a\quad (M_h+\sqrt{\alpha}K_h)v_b=\hat{r}_b.
  \end{equation*}
\end{itemize}
\end{algorithm}

Note that the matrix $G:=M_h+\sqrt{\alpha}K_h$ is symmetric positive definite. Hence, for the case where the (sparse) Cholesky factorizations of $G$ (need only to be done once) can be computed at a moderate cost, the above two linear system involving $G$ can be exactly and effectively solved. However, if the Cholesky factorizations of $G$ is not available, then the linear systems could be inexactly handled with some alternative efficient methods, e.g., preconditioned conjugate gradient (PCG) method, Chebyshev semi-iteration or some multigrid scheme. It is well known that the convergence behavior of iterative solution methods will
be precisely characterized in terms of $\kappa(M_h)$ and $\kappa(K_h)$, which represents the condition number of $M_h$ and $K_h$, respectively. Then about the bounds on the condition number, we have the following results, one can see Proposition 1.29 and Theorem 1.32 in \cite{spectral property} for more details.

\begin{theorem}\label{spectral property}
For $\mathcal{P}1$ approximation on a regular and quasi-uniform subdivision of $\mathbb{R}^n$ which satisfies Assumption {\rm\ref{regular and quasi-uniform triangulations}}, and for any $x\in \mathbb{R}^{N_h}$, the mass matrix $M_h$ approximates the scaled identity matrix in the sense that
\begin{equation*}
c_1 h^2\leq \frac{x^{T}M_hx}{x^{T}x}\leq c_2 h^2, \ if\  n=2, \ {\rm and}\ c_1 h^3\leq \frac{x^{T}M_hx}{x^{T}x}\leq c_2 h^3, \ if\  n=3.
\end{equation*}
The stiffness matrix $K_h$ satisfies
\begin{equation*}
d_1h^2\leq \frac{x^{T}K_hx}{x^{T}x}\leq d_2, \ if\  n=2, \ {\rm and }\ d_1h^3\leq \frac{x^{T}K_hx}{x^{T}x}\leq d_2 h, \ if\  n=3.
\end{equation*}
where the constants $c_1$, $c_2$, $d_1$ and $d_2$ are independent of the mesh size $h$.
\end{theorem}
Thus, according to Theorem \ref{spectral property}, in our numerical experiments, the approximation $\widehat{G}$ corresponding to the matrix $G:=M_h+\sqrt{\alpha}K_h$ is implemented by 20 steps of Chebyshev semi-iteration when the parameter $\alpha$ satisfies $\alpha\leq h^4$. Since in this case the coefficient matrix $G$ is dominated by the mass matrix and 20 steps of Chebyshev semi-iteration is an appropriate approximation for the action of $G$'s inverse. For the large values of $\alpha$, e.g., $\alpha>h^4$, however, the stiffness matrix $K_h$ makes a significant contribution. Hence, a fixed number of Chebyshev semi-iteration is no longer sufficient to approximate the action of $G^{-1}$. In this case, one typical choice is using a fixed number of algebraic multigrid ({\bf AMG}) V-cycles to approximate the action of $G^{-1}$. In our numerical implementation, the approximation $\widehat{G}$ to $G$ is set to be two {\bf AMG} V-cycles obtained by the amg operator in the iFEM software package\footnote{\noindent \textrm{For more details about the iFEM software package, we refer to the website \url{http://www.math.uci.edu/~chenlong/programming.html} }}.

In addition, let $(\hat {r}^k_1, \hat {r}^k_2)$ be the residual error vector, which means
\begin{equation}\label{equ:discretized saddle point problems3}
\left[
  \begin{array}{cc}
    \frac{1}{\alpha}M_h & \quad-K_h \\
    K_h & \quad M_h\\
  \end{array}
\right]\left[
         \begin{array}{c}
           \hat{p}^{k} \\
           \hat{y}^{k} \\
         \end{array}
       \right]=\left[
                 \begin{array}{c}
                  \frac{1}{\alpha}M_h(\lambda^k+\tilde{\mu}^k-\alpha y_r)+\hat {r}_1^k \\
                   M_hy_d+\hat {r}_2^k\\
                 \end{array}
               \right],
\end{equation}
and $\hat{\delta}_{p}^k=\hat {r}_1^k+K_hM_h^{-1}\hat {r}_2^k$. Thus in the numerical implementation we could require
\begin{equation}\label{error estimates12}
  \|\hat {r}^k_1\|_{2}+\|\hat {r}^k_2\|_{2}<\frac{\epsilon_k}{\max\{1,\|K_h\|_{2}\|M_h^{-1}\|_2\}},
\end{equation}
to guarantee the error vector $\|\hat{\delta}_{p}^k\|_{2} \leq {\epsilon_k}$.

\subsection{\textbf {An efficient predictor for the block $\tilde{p}$ subproblem}}\label{subsection5.3}
From the presentation in Step 1 of Algorithm \ref{algo1:Full inexact ABCD algorithm for (Dh)}, it appears that we would need to solve the block $p$ subproblem twice. In practice, in
order to improve the efficiency of our sGS-imABCD algorithm, in this section, we design an efficient predictor for the block $\tilde{p}$ subproblem to avoid solving it.

Obviously, to solve the block $\tilde{p}$ subproblem, we only need to replace $\lambda^k$ by $\tilde{\lambda}^k$ in the right-hand term of (\ref{equ:discretized saddle point problems2}). Then we have
\begin{equation}\label{equ:discretized saddle point problems4}
\left[\begin{array}{cc}
    \frac{1}{\alpha}M_h & \quad-K_h \\
    K_h & \quad M_h\\
  \end{array}
\right]\left[
         \begin{array}{c}
           \tilde{p}^{k} \\
           \tilde{y}^{k} \\
         \end{array}
       \right]=\left[
                 \begin{array}{c}
                  \frac{1}{\alpha}M_h(\tilde{\lambda}^k+\tilde{\mu}^k-\alpha y_r) \\
                   M_hy_d\\
                 \end{array}
               \right].
\end{equation}
Hence, all the numerical techniques for the block $\hat{p}$ is also applicable for the block $\tilde{p}$.

However, in practice, we can often avoid solving the linear system twice if $\hat{p}^k$ is already sufficiently close to $\tilde{p}^k$. More specifically, if we employ $\hat{p}^k$ to approximate $\tilde{p}^k$, then the residual vector for (\ref{equ:discretized saddle point problems4}) is given by
\begin{equation*}
  \left[\begin{array}{c}
  \tilde{r}_1^k \\
  \tilde{r}_2^k\\
  \end{array}
  \right]=\left[\begin{array}{c}
  \frac{1}{\alpha}M_h(\tilde{\lambda}^k-\lambda^k)-\hat{r}_1^k \\
  -\hat{r}_2^k\\
  \end{array}
  \right],
\end{equation*}
which means $\tilde{\delta}_{p}^k=\frac{1}{\alpha}M_h(\tilde{\lambda}^k-\lambda^k)-\hat{r}_1^k-K_hM_h^{-1}\hat {r}_2^k$. If the condition
\begin{equation}\label{error estimates2}
  \|\tilde {r}^k_1\|_{2}+\|\tilde {r}^k_2\|_{2}<\frac{\epsilon_k}{\max\{1,\|K_h\|_{2}\|M_h^{-1}\|_2\}},
\end{equation}
is satisfied which can also guarantee the error vector $\|\tilde{\delta}_{p}^k\|_{2} \leq {\epsilon_k}$, then we need not solve the linear system (\ref{equ:discretized saddle point problems4}) and take $\tilde {p}^k=\hat{p}^k$.

At last, although we solve problem (\ref{eqn:orginal problems}) via its dual, our ultimate goal is look for optimal control solution. Thus we should introduce the KKT condition for (\ref{eqn:discretized matrix-vector dual problem}) as below
\begin{equation*}
\left\{
\begin{array}{r@{\;=\;}l}
0 & M_h(y- y_d)+ K_hp,\\
0 & \alpha u -p + \lambda + \mu,\\
0 & K_hy- M_hu-M_h y_r,\\
0 & u- {\rm\Pi}_{[a,b]}(u+M_h\mu),\\
0 & \lambda- {\rm\Pi}_{[-\beta,\beta]}(\lambda+ M_hu).
\end{array}
\right.
\end{equation*}
Then we can have $u= (p - \lambda - \mu)/\alpha$.

Furthermore, in order to measure the accuracy of an approximate optimal solution $(\mu,\lambda,p)$ for (\ref{eqn:discretized matrix-vector dual problem}), let us introduce the checkable stopping criterion for our sGS-ABCD algorithm. Let $\epsilon$ be a given accuracy tolerance, we terminate our sGS-imABCD method when $\eta<\epsilon$,
where the relative residual $\eta$ is given by
\begin{equation}\label{residual}
  \eta=\max{\{\eta_1,\eta_2,\eta_3,\eta_4\}},
\end{equation}
where
\begin{equation*}
  \begin{aligned}
    \eta_1=\frac{\|M_h(y- y_d)+ K_hp\|}{1+\|M_hy_d\|}, & \quad \eta_2=\frac{\|K_hy- M_hu-M_h y_r\|}{1+\|M_hy_r\|},\\
    \eta_3=\frac{\|u- {\rm\Pi}_{[a,b]}(u+M_h\mu)\|}{1+\|u\|},& \quad\eta_4=\frac{\|\lambda- {\rm\Pi}_{[-\beta,\beta]}(\lambda+ M_hu)\|}{1+\|\lambda\|},
  \end{aligned}
\end{equation*}
and $u= (p - \lambda - \mu)/\alpha$.

\section{\textbf {An ihADMM method and an APG method for (\ref{equ:approx discretized matrix-vector form})}}\label{subsection5.4}
In this section, we will introduce some algorithms for comparison. First, as already mentioned, in order to show the efficiency of the duality-based approach to solve problem (\ref{eqn:discretized matrix-vector dual problem}), we also use our sGS-imABCD method to solve problem ({\ref{eqn:approx discretized matrix-vector dual problem}}) for comparison.

Comparing ({\ref{eqn:approx discretized matrix-vector dual problem}}) with (\ref{eqn:discretized matrix-vector dual problem}), we can easily see that our sGS-imABCD method applied to ({\ref{eqn:approx discretized matrix-vector dual problem}}) is almost the same as that for (\ref{eqn:discretized matrix-vector dual problem}), except the $\lambda$-subproblem. For the  $\lambda$-subproblem, we have
\begin{equation*}
  \tilde{\lambda}^{k}
        =\arg\min\frac{1}{2\alpha}\|\lambda-(\hat{p}^{k}-\tilde{\mu}^k)\|_{M_h}^2+\delta_{[-\beta,\beta]}(W_h^{-1}M_h\lambda)+\frac{1}{2\alpha}\|\lambda- \lambda^{k}\|_{c_n M_hW_h^{-1}M_h-M_h}^2.
\end{equation*}
Let $d=M_h\lambda$, then we have,
\begin{eqnarray*}
  \tilde{d}^{k}
        &=&\arg\min\frac{1}{2\alpha}\|d-M_h(\hat{p}^{k}-\tilde{\mu}^k)\|_{M_h^{-1}}^2+\delta_{[-\beta,\beta]}(W_h^{-1}d)+\frac{1}{2\alpha}\|d- d^{k}\|_{c_n W_h^{-1}-M_h^{-1}}^2,\\
        &=&W_h{\rm\Pi}_{[-\beta,\beta]}(e^k),
\end{eqnarray*}
where $e^k:=W_h^{-1}M_h\lambda^k+\frac{1}{c_n}(\hat {p}^k-\tilde{\mu}^k-\lambda^k)$. Then we have
$\tilde\lambda^k=M_h^{-1}\tilde{d}^k$.
\begin{remark}
Similar we could obtain the KKT equation for problem (\ref{equ:approx discretized matrix-vector form}) and ({\ref{eqn:approx discretized matrix-vector dual problem}}) as below
\begin{equation*}
\left\{
\begin{array}{r@{\;=\;}l}
0 & M_h(y- y_d)+ K_h p,\\
0 & \alpha u -p + \lambda+\mu,\\
0 & K_h y- M_h u-M_h y_r,\\
0 & u- {\rm\Pi}_{U_{ad}}(u+M_h\mu),\\
0 & M_h\lambda- W_h{\rm\Pi}_{[-\beta, \beta]}(W_h^{-1}M_h\lambda+ W_h u).
\end{array}
\right.
\end{equation*}
Thus we measure the accuracy of an approximate optimal solution $(\mu,\lambda,p)$ for (\ref{eqn:discretized matrix-vector dual problem}) by using the following relative residual:
\begin{equation}\label{residual2}
  \eta=\max{\{\eta_1,\eta_2,\eta_3,\eta_4\}},
\end{equation}
where
\begin{equation*}
  \begin{aligned}
    \eta_1=\frac{\|M_h(y- y_d)+ K_hp\|}{1+\|M_hy_d\|},& \quad \eta_2=\frac{\|K_hy- M_hu-M_h y_r\|}{1+\|M_hy_r\|}, \\
    \eta_3=\frac{\|u- {\rm\Pi}_{[a,b]}(u+ M_h\mu)\|}{1+\|u\|}, &\quad \eta_4=\frac{\|M_h\lambda- W_h{\rm\Pi}_{[-\beta, \beta]}(W_h^{-1}M_h\lambda+  W_h u)\|}{1+\|M_h\lambda\|},
   \end{aligned}
\end{equation*}
and $u=(p - \lambda - \mu)/\alpha$. We also terminate our sGS-imABCD method for ({\ref{eqn:approx discretized matrix-vector dual problem}}) when $\eta<\epsilon$.
\end{remark}

Instead of the sGS-imABCD method, one can also apply the ihADMM \cite{iwADMM} and APG method \cite{FIP} to solve the primal problem (\ref{equ:approx discretized matrix-vector form}) for the sake of numerical comparison, and the details are given as follows.
\begin{algorithm}[H]\label{algo4:inexact weighted ADMM for problem RHP}
  \caption{\textbf {inexact heterogeneous ADMM (ihADMM) algorithm for (\ref{equ:approx discretized matrix-vector form})}}
  \KwIn{$(z^0, u^0, \lambda^0)\in {\rm dom} (\delta_{[a,b]}(\cdot))\times \mathbb{R}^n \times \mathbb{R}^n $ and a parameter $\tau \in (0,1]$. Let $\{\epsilon_k\}^\infty_{k=0}$ be a sequence satisfying $\{\epsilon_k\}^\infty_{k=0}\subseteq [0,+\infty)$ and $\sum\limits_{k=0}^{\infty}\epsilon_k<\infty$. Set $k=0$}
  \KwOut{$ u^k, z^{k}, \lambda^k$}
\begin{description}
\item[\bf Step 1] Solving the follwing linear system (inexact)
\begin{equation*}
\left[
  \begin{array}{cc}
    (\alpha+\sigma)K_h & M_h \\
    -M_h & K_h\\
  \end{array}
\right]\left[
         \begin{array}{c}
           u^{k+1} \\
           y^{k+1} \\
         \end{array}
       \right]\approx\left[
                 \begin{array}{c}
                   M_hy_d+\sigma K_hz^k-K_h\lambda^k \\
                   M_h y_r\\
                 \end{array}
               \right],
\end{equation*}
with the residual error vector ${r}^k=(r^k_1;r^k_2)$ satisfies $\|{r}^k\|_{2} \leq {\epsilon_k}$
\item[\bf Step 2] Compute $z^k$ as follows:
       \begin{eqnarray*}
       z^{k+1}&=&{\rm\Pi}_{[a,b]}({\rm soft}(u^{k+1}+(W_h^{-1}M_h\lambda^k)/\sigma), \beta/\sigma),
       \end{eqnarray*}
  \item[\bf Step 3] Compute
  \begin{eqnarray*}
    \lambda^{k+1} &=& \lambda^k+\tau\sigma(u^{k+1}-z^{k+1}).
  \end{eqnarray*}

  \item[\bf Step 4] If a termination criterion is not met, set $k:=k+1$ and go to Step 1
\end{description}
\end{algorithm}
\begin{remark}\label{remark3}
To solve the linear system in Step 1 of Algorithm \ref{algo4:inexact weighted ADMM for problem RHP}, we also use the GMRES method with the {\rm PMHSS} preconditioner. Thus all the numerical techniques as mentioned in Subsection \ref{subsection5.2} can be used. 
\end{remark}

\begin{algorithm}[H]\label{algo1:APG algorithm for (Ph)}
  \caption{\textbf{(APG algorithm for (\ref{equ:approx discretized matrix-vector form}))}}
 \KwIn{$\tilde u^1=u^0\in {\rm dom} (\delta_{[a,b]})$, $\rho>1$ and $L^0>0$. Let $\{\epsilon_k\}^\infty_{k=0}$ be a sequence satisfying $\{\epsilon_k\}^\infty_{k=0}\subseteq [0,+\infty)$ and  $\sum\limits_{k=0}^{\infty}\epsilon_k<\infty$. Set $k= 1, t_1= 1.$}
 \KwOut{$ (y^k, u^k, p^k)$}
 Iterate until convergence
\begin{description}
  \item[\bf Step 1] Choose error tolerance $\tilde\delta_{y}^k, \tilde\delta_{p}^k, \delta_y^k,\delta_p^k$ such that
  \begin{equation*}
\max\{\|{\tilde\delta}_{y}^k|\|,\|\tilde{\delta}_p^k|\|,\|\delta_y^k|\|,\|\delta_p^k|\|\}\leq \epsilon_k.
  \end{equation*}
\begin{itemize}
\item[]Compute 
\begin{equation*}
  K_h\tilde {y}^k\approx M_h\tilde {u}^k+M_h y_r,\quad K_h\tilde {p}^k\approx M_h y_d-M_h \tilde {y}^k,
\end{equation*}
with the residual error vector $\tilde\delta_{y}^k$ and $\tilde\delta_{p}^k$, respectively.
  \item[\bf Step 2] Backtracking: Find the smallest nonnegative integer $i$ such that with $L=\rho^iL^{k-1}$
           \begin{equation*}
             \hat{J}(v)\leq \hat{J}(\tilde {u}^k)+\langle\alpha M_h\tilde {u}^k-M_h \tilde {p}^k,  v-\tilde {u}^k\rangle+\frac{L}{2}\|v-\tilde {u}^k\|^2,
           \end{equation*}
where
\begin{eqnarray*}
&&\hat{J}(v)=\frac{1}{2}\|y-y_{d}\|_{M_h}^{2}+\frac{\alpha}{2}\|v\|_{M_h}^{2},\\
&&v={\rm\Pi}_{[a,b]}(W_h{\rm soft}(W_h^{-1}(\tilde {u}^{k}-(\alpha M_h\tilde{u}^k- M_h\tilde {p}^k)/L), L/\beta)),\\
&&K_h{y}\approx M_hv+M_h y_r.
\end{eqnarray*}
 \item[\bf Step 3] Set $L^k=L$ and compute
\begin{eqnarray*}
u^k={\rm\Pi}_{[a,b]}(W_h{\rm soft}(W_h^{-1}(\tilde {u}^{k}-(\alpha M_h\tilde{u}^k-\tilde M_hp^k)/L^k), \beta/L^k),
\end{eqnarray*}
\item[]and 
\begin{equation*}
  K_h{y}^k\approx M_h{u}^k+M_h y_r,\quad  K_h{p}^k\approx M_h y_d-M_h {y}^k,
\end{equation*}
with the residual error vector $\delta_{y}^k$ and $\delta_{p}^k$, respectively.
\end{itemize}
  \item[\bf Step 4] Set $t_{k+1}=\frac{1+\sqrt{1+4t_k^2}}{2}$ and $\beta_k=\frac{t_k-1}{t_{k+1}}$, Compute
\begin{eqnarray*}
\tilde {u}^{k+1}=u^{k}+ \beta_{k}(u^{k}-u^{k-1}).
\end{eqnarray*}
\end{description}
\end{algorithm}

\begin{remark}\label{remark4}
To inexactly solve the linear system about the coefficient matrix $K_h$, in our implementations, we use two AMG V-cycles method to approximate $K_h$.
\end{remark}

\begin{remark}\label{remark5}
Moreover, for problem (\ref{equ:approx discretized matrix-vector form}), the the KKT condition can be given by
\begin{equation*}
\left\{
\begin{array}{r@{\;=\;}l}
0 & M_h(y- y_d)+ K_h p,\\
0 & \alpha u -p + \lambda,\\
0 & K_h y- M_h u-M_h y_r,\\
0 & u-{\rm\Pi}_{[a,b]}(W_h{\rm soft}(W_h^{-1}({u}+M_h\lambda), \beta).
\end{array}
\right.
\end{equation*}
Let $\epsilon$ be a given accuracy tolerance. Thus we terminate both ihADMM method and APG method when $\eta<\epsilon$,
where
\begin{equation}\label{residual3}
  \eta=\max{\{\eta_1,\eta_2,\eta_3\}},
\end{equation}
in which
\begin{equation*}
     \eta_1=\frac{\|M_h(y- y_d)+ K_hp\|}{1+\|M_hy_d\|},\
     \eta_2=\frac{\|K_hy- M_hu-M_h y_r\|}{1+\|M_hy_r\|},\
     \eta_3=\frac{\|u-{\rm\Pi}_{[a,b]}(W_h{\rm soft}(W_h^{-1}({u}+M_h\lambda), \beta)\|}{1+\|u\|},
\end{equation*}
$p=\alpha u + \lambda$ for ihADMM method and $\lambda=p-\alpha u$ for APG method.
\end{remark}

\section{Numerical results}\label{sec:6}
In this section, we will use the following examples to evaluate the numerical behaviour of our sGS-imABCD method for (\ref{eqn:discretized matrix-vector dual problem}) and verify the theoretical error estimates given in Theorem \ref{theorem:error1}. For comparison, we will also show the numerical results obtained by the our sGS-imABCD for ({\ref{eqn:approx discretized matrix-vector dual problem}}) and the ihADMM and APG methods for (\ref{equ:approx discretized matrix-vector form}).

\subsection{\textbf {Algorithmic details}}\label{subsec:6.1}
We begin by describing the algorithmic details which are common to all examples, unless otherwise mentioned.

\textbf{Discretization.} As show in Section \ref{sec:3}, the discretization was carried out using piece-wise linear and continuous finite elements.
The assembly of mass and the stiffness matrices, as well as the lump mass matrix was left to the iFEM software package.

To present the finite element error estimates results, it is convenient to introduce the experimental order of convergence (EOC), which for some positive error functional $E(h)$ with $h> 0$ is defined as follows: Given two grid sizes $h_1\neq h_2$, let
\begin{equation}\label{EOC}
  \mathrm{EOC}:=\frac{\log {E(h_1)}-\log {E(h_2)}}{\log{h_1}-\log {h_2}}.
\end{equation}
It follows from this definition that if $E(h)=\mathcal{O}(h^{\gamma})$ then $\mathrm{EOC}\approx\gamma$. The error functional $E(\cdot)$ investigated in the present section is given by
\begin{equation}\label{error norm}
  E_2(h):=\|u-u_h\|_{L^2{(\Omega)}}.
\end{equation}

\textbf{Initialization.} For all numerical examples, we choose the initial values as zero for all algorithms.

\textbf{Parameter setting.} For the ihADMM method, the step-length $\tau$ for lagrangian multipliers $\lambda$ was chosen as $\tau=1$, and the penalty parameter $\sigma$ was chosen as $\sigma=0.1 \alpha$. For the APG method, we estimate an approximation to the Lipschitz constant $L$ with a backtracking method with $\eta=1.4$ and $L^0=10^{-8}$.

\textbf{Stopping criterion.} In our numerical experiments, we terminate all the algorithms when the corresponding relative residual $\eta<10^{-7}$.

\textbf{Computational environment.}
All our computational results are obtained by MATLAB Version 8.5(R2015a) running on a computer with
64-bit Windows 7.0 operation system, Intel(R) Core(TM) i7-5500U CPU (2.40GHz) and 8GB of memory.

\subsection{\textbf {Examples}}\label{subsec:6.2}
Before giving the specific examples, we first introduce the following procedure, which can help us formulate sparse optimal control problems.

\begin{algorithm}[H]\label{algo:construct the optimal problem}
\caption{\textbf{Construct the optimal control problem}}
  \begin{description}
  \item [\bf Step 1]. Choose $y^*\in H^1_0(\Omega)$ and $p^*\in H^1_0(\Omega)$ arbitrarily.
  \item [\bf Step 2]. Set
  \begin{equation*}
    u^*:=\left\{\begin{aligned}
       &\min\{\frac{p^*-\beta}{\alpha}, b\}, &on~~{x\in\Omega: p^*(x)>\beta},  \\
       &\max\{\frac{p^*+\beta}{\alpha}, a\}, \quad&on~~{x\in\Omega: p^*(x)<-\beta},\\
       &0,                                      &elsewhere.
    \end{aligned}\right.
  \end{equation*}
  \item [\bf Step 3]. Set $y_r=Ay^*-Bu^*$ and $y_d=Ap^*+y^*$.
\end{description}
\end{algorithm}
According to the first-order optimality condition in Theorem \ref{First-order optimality condition}, we can see that Algorithm \ref{algo:construct the optimal problem} provides an optimal solution $(y^*, u^*)$ of the sparse optimal control problem (\ref{eqn:orginal problems}). Thus we can construct examples for which we know the exact solution through the above procedure.
\begin{example}\label{example:1}
  \begin{equation*}
     \left\{ \begin{aligned}
        &\min \limits_{(y,u)\in H^1_0(\Omega)\times L^2(\Omega)}^{}\ \ J(y,u)=\frac{1}{2}\|y-y_d\|_{L^2(\Omega)}^{2}+\frac{\alpha}{2}\|u\|_{L^2(\Omega)}^{2}+\beta\|u\|_{L^1(\Omega)} \\
        &\qquad\quad{\rm s.t.}\qquad \quad\quad\quad-\Delta y=u+y_r\quad \mathrm{in}\  \Omega,\\
         &\qquad \qquad \quad\qquad \qquad \qquad  ~y=0\quad  \mathrm{on}\ \partial\Omega,\\
         &\quad\qquad \qquad \qquad\qquad  u\in U_{ad}=\{v(x)|a\leq v(x)\leq b, {\rm a.e }\  \mathrm{on}\ \Omega \}.
                          \end{aligned} \right.
 \end{equation*}
Here, we consider the problem with control $u\in L^2(\Omega)$ on the unit square $\Omega= (0, 1)^2$ with $\alpha=0.5$, $\beta=0.5$, $a=-0.5$ and $b=0.5$. It is a constructed problem, thus we set $y^*=\sin(2\pi x_1)\exp(0.5x_1)\sin(4\pi x_2)$ and $p^*=2\beta\sin(2\pi x_1)\exp(0.5x_1)\sin(4\pi x_2)$. Then through Algorithm \ref{algo:construct the optimal problem}, we can easily get the optimal control solution $u^*$, the source term $y_r$ and the desired state $y_d$.

An example for the discretized optimal control on mesh $h=2^{-7}$ is shown in Figure \ref{example1fig:control on h=$2^{-7}$}. The error of the control $u$ w.r.t the $L^2$ norm and the experimental order of convergence (EOC) for control are presented in Table \ref{tab:1} and Table \ref{tab:3}. They also confirm that indeed the convergence rate is of order $O(h)$. Compared the error results from Table \ref{tab:1} and Table \ref{tab:3}, it is obvious to see that solving the dual problem (\ref{eqn:discretized matrix-vector dual problem}) could get better error results than that from solving (\ref{equ:approx discretized matrix-vector form}) and ({\ref{eqn:approx discretized matrix-vector dual problem}}).

Numerical results for the accuracy of solution, number of iterations and cpu time obtained by the our proposed sGS-imABCD method for (\ref{eqn:discretized matrix-vector dual problem}) are also shown in Table \ref{tab:1}. As a result we obtain from Table \ref{tab:1}, one can see that our proposed sGS-imABCD method is an efficient algorithm to solve problem (\ref{eqn:discretized matrix-vector dual problem}) to high accuracy. It should be pointed out that iter.$\tilde{p}$-block denotes the iterations of $\tilde{p}$ in Table \ref{tab:1}. It is clear that $p$-subproblem almost always not be computed twice, which demonstrates the efficiency of our strategy to predict the solution of $\tilde{p}$-subproblem. Furthermore, the numerical results in terms of iteration numbers illustrate the mesh-independent performance of our proposed sGS-imABCD method. Additionally, in Table \ref{tab:2}, we list the numbers of iteration steps and the relative residual errors of PMHSS-preconditioned GMRES method for the $\hat{p}$-subproblem on mesh $h=2^{-7}$ and $h=2^{-8}$. From Table \ref{tab:2}, we can see that the number of iteration steps of the PMHSS-preconditioned GMRES method is roughly independent of the mesh size $h$.

As a comparison, numerical results obtained by the our proposed sGS-imABCD method for ({\ref{eqn:approx discretized matrix-vector dual problem}}) and the iwADMM and APG methods for (\ref{equ:approx discretized matrix-vector form}) are shown in Table \ref{tab:3}. As a result from Table \ref{tab:3}, it can be observed that our sGS-imABCD is faster and more efficient than the iwADMM and APG methods in terms of the iterations and CPU times.

At last, in order to show the robustness of our proposed sGS-imABCD method with respect to the parameters $\alpha$ and $\beta$, we also test the same problem with different values of $\alpha$ and $\beta$ on mesh $h=2^{-8}$. The results are presented in Table \ref{tab:4}. From the Table \ref{tab:4}, it is obviouse to see that our method could solve problem (\ref{eqn:discretized matrix-vector dual problem}) to high accuracy for all tested values of $\alpha$ and $\beta$ within 50 iterations. More importantly, from the results, we can see that when $\alpha$ is fixed, the number of iteration steps of the sGS-imABCD method remains nearly constant for $\beta$ ranging from $0.005$ to $1$. However, for a fixed $\beta$, as $\alpha$ increases from $0.005$
to $0.5$, the number of iteration steps of the sGS-imABCD method changes drastically. These observations indicate that the sGS-imABCD method shows the $\beta$-independent convergence property, whereas it dose not have the same convergence property with respect to the parameter $\alpha$. It should be pointed out that the numerical results are also consistent with the theoretical conclusion which based on Theorem \ref{sGS-imABCD convergence}.

\begin{figure}[H]
\begin{center}
\includegraphics[width=0.40\textwidth]{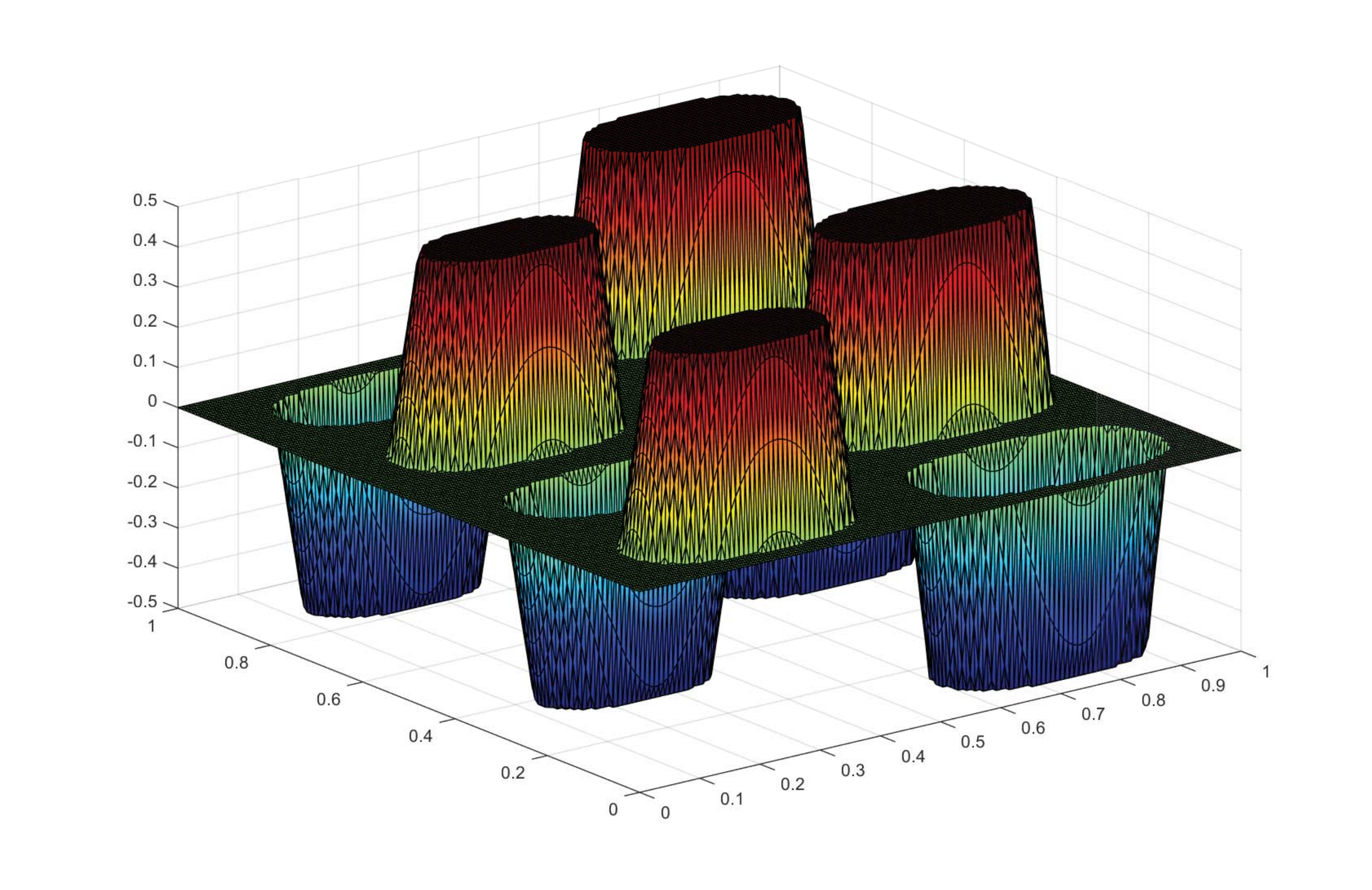}
\includegraphics[width=0.40\textwidth]{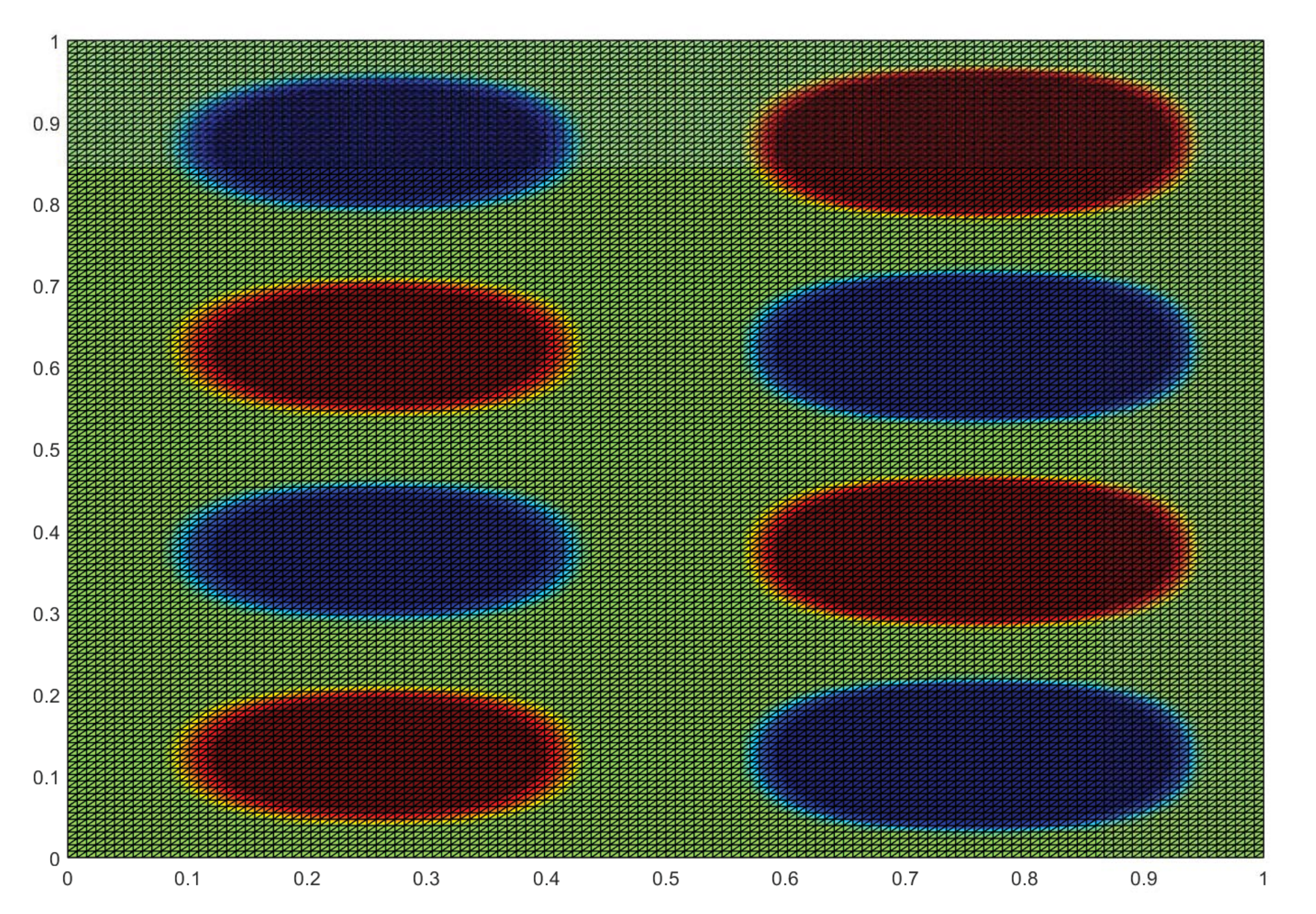}
\caption{Optimal control $u_h$ on the square, $h=2^{-7}$. Dark red and dark blue areas correspond to $u_h=\pm 0.5$ and green areas to $u_h=0$}\label{example1fig:control on h=$2^{-7}$}.
\end{center}
\end{figure}

\begin{table}[H]\footnotesize
\caption{Example \ref{example:1}: The performance of sGS-imABCD for (\ref{eqn:discretized matrix-vector dual problem}). In the table, $\#$dofs stands for the number of degrees of freedom for the control variable on each grid level.}\label{tab:1}
\begin{center}
\begin{tabular}{@{\extracolsep{\fill}}ccccccccccc}
\hline
\multirow{2}{*}{$h$}&\multirow{2}{*}{$\#$dofs} &\multirow{2}{*}{iter.sGS-imABCD}&\multirow{2}{*}{iter.$\tilde{p}$-block}&\multirow{2}{*}{}& \multirow{2}{*}{residual $\eta$}&\multirow{2}{*}{} & \multirow{2}{*}{CPU time/s} & \multirow{2}{*}{$E_2$}&\multirow{2}{*}{EOC}\\
& &  & & &  & &  &  &\\
\hline
&&&&&&&\\
$2^{-3}$  &49 &13 &4  &   &6.60e-08  &   &0.14 &0.1784      &--       &    \\
&&&&&&&\\
$2^{-4}$  &225 &13 &4  &   &6.32e-08  &   &0.20 &0.0967      &0.8834       &    \\
&&&&&&&\\
$2^{-5}$  &961 &12 &3  &   &7.38e-08  &   &0.33 &0.0399      &1.0803   &    \\
&&&&&&&\\
$2^{-6}$  &3969 &13 &3  &   &9.78e-08  &   &2.04 &0.0155    &1.1749   &    \\
&&&&&&&\\
$2^{-7}$  &16129&12 &3  &   &6.66e-08  &   &8.25 &0.0052    &1.2754  &    \\
&&&&&&&\\
$2^{-8}$  &65025&10 &3  &   &7.05e-08  &   &52.15 &0.0017   &1.3388   &    \\
&&&&&&&\\
$2^{-9}$  &261121&9 &2 &   &5.19e-08  &   &312.82 &0.0006  & 1.3617   &    \\
&&&&&&&\\
\hline
\end{tabular}
\end{center}
\end{table}

\begin{table}[H]\footnotesize
\caption{Example \ref{example:1}: The convergence behavior of GMRES for $\hat{p}$-block subproblem. 
}\label{tab:2}
\begin{center}
\begin{tabular}{@{\extracolsep{\fill}}|c|c|c|c|}
\hline
\multirow{2}{*}{$h$}&\multirow{2}{*}{iter.sGS-imABCD}&\multirow{2}{*}{iter.GMRES of $\hat{p}$-block}&\multirow{2}{*}{Relative residual error of GMRES} \\
&&&\\
\hline

                            &\multirow{2}{*}{1}    &\multirow{2}{*}{8} &\multirow{2}{*}{1.30e-07} \\
                            &\multirow{2}{*}{2}    &\multirow{2}{*}{4} &\multirow{2}{*}{1.07e-07}\\
                            &\multirow{2}{*}{3}    &\multirow{2}{*}{4} &\multirow{2}{*}{5.26e-08}\\
                            &\multirow{2}{*}{4}    &\multirow{2}{*}{4} &\multirow{2}{*}{1.56e-08}\\
                            &\multirow{2}{*}{5}    &\multirow{2}{*}{4} &\multirow{2}{*}{2.05e-09}\\
                            &\multirow{2}{*}{6}    &\multirow{2}{*}{4} &\multirow{2}{*}{1.58e-09}\\
\multirow{2}{*}{$2^{-7}$}   &\multirow{2}{*}{7}    &\multirow{2}{*}{4} &\multirow{2}{*}{1.23e-09}\\
                            &\multirow{2}{*}{8}    &\multirow{2}{*}{4} &\multirow{2}{*}{1.29e-10}\\
                            &\multirow{2}{*}{9}    &\multirow{2}{*}{2} &\multirow{2}{*}{1.16e-10}\\
                            &\multirow{2}{*}{10}   &\multirow{2}{*}{2} &\multirow{2}{*}{1.07e-10}\\
                            &\multirow{2}{*}{11}   &\multirow{2}{*}{2} &\multirow{2}{*}{5.98e-11}\\
                            &\multirow{2}{*}{12}   &\multirow{2}{*}{2} &\multirow{2}{*}{1.30e-11}\\
                            &&&\\
\hline

                            &\multirow{2}{*}{1}    &\multirow{2}{*}{8} &\multirow{2}{*}{6.31e-08}\\
                            &\multirow{2}{*}{2}    &\multirow{2}{*}{4} &\multirow{2}{*}{2.18e-08}\\
                            &\multirow{2}{*}{3}    &\multirow{2}{*}{4} &\multirow{2}{*}{8.43e-09}\\
                            &\multirow{2}{*}{4}    &\multirow{2}{*}{4} &\multirow{2}{*}{3.18e-09}\\
                            &\multirow{2}{*}{5}    &\multirow{2}{*}{4} &\multirow{2}{*}{1.07e-09}\\
\multirow{2}{*}{$2^{-8}$}   &\multirow{2}{*}{6}    &\multirow{2}{*}{4} &\multirow{2}{*}{5.53e-10}\\
                            &\multirow{2}{*}{7}    &\multirow{2}{*}{4} &\multirow{2}{*}{5.25e-11}\\
                            &\multirow{2}{*}{8}    &\multirow{2}{*}{4} &\multirow{2}{*}{5.90e-12}\\
                            &\multirow{2}{*}{9}    &\multirow{2}{*}{2} &\multirow{2}{*}{4.86e-12}\\
                            &\multirow{2}{*}{10}    &\multirow{2}{*}{2} &\multirow{2}{*}{4.18e-12}\\
                            &&&\\
\hline
\end{tabular}
\end{center}
\end{table}

\begin{table}[H]\footnotesize
\caption{Example \ref{example:1}: The convergence behavior of sGS-imABCD for ({\ref{eqn:approx discretized matrix-vector dual problem}}) , ihADMM and APG for (\ref{equ:approx discretized matrix-vector form}). In the table, $\#$dofs stands for the number of degrees of freedom for the control variable on each grid level. $E_2=\min\{E_2(sGS-imABCD), E_2(ihADMM), E_2(APG)\}$. 
}\label{tab:3}
\begin{center}
\begin{tabular}{@{\extracolsep{\fill}}|c|c|c|c|c|c|c|c|}
\hline
\multirow{2}{*}{$h$}&\multirow{2}{*}{$\#$dofs}&\multirow{2}{*}{$E_2$}&\multirow{2}{*}{EOC}& \multirow{2}{*}{Index of performance} &\multirow{2}{*}{sGS-imABCD} & \multirow{2}{*}{ihADMM} & \multirow{2}{*}{APG} \\
&&&&&&&\\
\hline

                            &&&&\multirow{2}{*}{iter}            &\multirow{2}{*}{13} &\multirow{2}{*}{32} &\multirow{2}{*}{16}              \\
\multirow{2}{*}{$2^{-3}$}   &\multirow{2}{*}{49}&\multirow{2}{*}{0.2925}&\multirow{2}{*}{--}
                            &\multirow{2}{*}{residual $\eta$}
                            &\multirow{2}{*}{6.25e-08}            &\multirow{2}{*}{6.33e-08}
                            &\multirow{2}{*}{3.51e-08}             \\
                            &&&&\multirow{2}{*}{CPU time/s}          &\multirow{2}{*}{0.16} &\multirow{2}{*}{0.23} &\multirow{2}{*}{0.22}             \\
                            &&&&&&&\\
\hline
                            &&&&\multirow{2}{*}{iter}            &\multirow{2}{*}{12} &\multirow{2}{*}{36} &\multirow{2}{*}{18}             \\
\multirow{2}{*}{$2^{-4}$}   &\multirow{2}{*}{225}&\multirow{2}{*}{0.1127}&\multirow{2}{*}{1.3759}&\multirow{2}{*}{residual $\eta$}
                            &\multirow{2}{*}{6.34e-08}            &\multirow{2}{*}{8.91e-08}
                            &\multirow{2}{*}{7.23e-08}             \\
                            &&&&\multirow{2}{*}{CPU times/s}             &\multirow{2}{*}{0.24} &\multirow{2}{*}{0.44} &\multirow{2}{*}{0.45}  \\
                            &&&&&&&\\
\hline
                            &&&&\multirow{2}{*}{iter}            &\multirow{2}{*}{13} &\multirow{2}{*}{40} &\multirow{2}{*}{16}            \\
\multirow{2}{*}{$2^{-5}$}   &\multirow{2}{*}{961}&\multirow{2}{*}{0.0457}&\multirow{2}{*}{1.3390}
                            &\multirow{2}{*}{residual $\eta$}
                            &\multirow{2}{*}{7.10e-08}&\multirow{2}{*}{7.42e-08}
                            &\multirow{2}{*}{8.88e-08}             \\
                            &&&&\multirow{2}{*}{CPU time/s}          &\multirow{2}{*}{0.47} &\multirow{2}{*}{1.17} &\multirow{2}{*}{2.98}            \\
                            &&&&&&&\\
\hline
                            &&&&\multirow{2}{*}{iter}            &\multirow{2}{*}{14} &\multirow{2}{*}{44} &\multirow{2}{*}{16}              \\
\multirow{2}{*}{$2^{-6}$}   &\multirow{2}{*}{3969}&\multirow{2}{*}{0.0161}&\multirow{2}{*}{1.3944}&\multirow{2}{*}{residual $\eta$}
                            &\multirow{2}{*}{4.05e-08}&\multirow{2}{*}{9.10e-08}
                            &\multirow{2}{*}{6.60e-08}             \\
                            &&&&\multirow{2}{*}{CPU time/s}          &\multirow{2}{*}{2.62} &\multirow{2}{*}{6.04} &\multirow{2}{*}{4.86}            \\
                            &&&&&&&\\
\hline
                            &&&&\multirow{2}{*}{iter}            &\multirow{2}{*}{12} &\multirow{2}{*}{50} &\multirow{2}{*}{16}              \\
\multirow{2}{*}{$2^{-7}$}   &\multirow{2}{*}{16129}&\multirow{2}{*}{0.0058}&\multirow{2}{*}{1.4132}&\multirow{2}{*}{residual $\eta$}
                            &\multirow{2}{*}{ 6.43e-08}&\multirow{2}{*}{9.80e-08}
                            &\multirow{2}{*}{8.45e-08}             \\
                            &&&&\multirow{2}{*}{CPU time/s}          &\multirow{2}{*}{10.22} &\multirow{2}{*}{29.53} &\multirow{2}{*}{30.63}            \\
                            &&&&&&&\\
\hline
                            &&&&\multirow{2}{*}{iter}            &\multirow{2}{*}{10} &\multirow{2}{*}{53} &\multirow{2}{*}{17}              \\
\multirow{2}{*}{$2^{-8}$}   &\multirow{2}{*}{65025}&\multirow{2}{*}{0.0019}&\multirow{2}{*}{1.4503}
                            &\multirow{2}{*}{residual $\eta$}
                            &\multirow{2}{*}{7.05e-08}&\multirow{2}{*}{8.93e-08}
                            &\multirow{2}{*}{8.88e-08}             \\
                            &&&&\multirow{2}{*}{CPU time/s}          &\multirow{2}{*}{60.45} &\multirow{2}{*}{160.24} &\multirow{2}{*}{92.60}            \\
                            &&&&&&&\\
\hline
                            &&&&\multirow{2}{*}{iter}            &\multirow{2}{*}{10} &\multirow{2}{*}{54} &\multirow{2}{*}{18}              \\
\multirow{2}{*}{$2^{-9}$}   &\multirow{2}{*}{261121}&\multirow{2}{*}{0.0007}&\multirow{2}{*}{1.4542}
                            &\multirow{2}{*}{residual $\eta$}
                            &\multirow{2}{*}{5.21e-08}&\multirow{2}{*}{7.96e-08}
                            &\multirow{2}{*}{3.24e-08}             \\
                            &&&&\multirow{2}{*}{CPU time/s}          &\multirow{2}{*}{395.78} &\multirow{2}{*}{915.71} &\multirow{2}{*}{859.22}            \\
                            &&&&&&&\\
\hline
\end{tabular}
\end{center}
\end{table}

\begin{table}[H]\footnotesize
\caption{Example \ref{example:1}: The performance of sGS-imABCD for (\ref{eqn:discretized matrix-vector dual problem}) with different values of $\alpha$ and $\beta$.}\label{tab:4}
\begin{center}
\begin{tabular}{@{\extracolsep{\fill}}ccccccccccc}
\hline
&\multirow{2}{*}{$h$}&\multirow{2}{*}{}&\multirow{2}{*}{$\alpha$}&\multirow{2}{*}{}&\multirow{2}{*}{$\beta$}&\multirow{2}{*}{} &\multirow{2}{*}{iter.sGS-imABCD}&\multirow{2}{*}{}&\multirow{2}{*}{residual error $\eta$ about K-K-T} \\
&&&& & &  & & &\\
\hline
&&&& & &  & & &\\
  &&&&&0.005& &49& &7.59e-08      \\
&&&0.005&  &0.05& &48& &8.86e-08    \\
  &&&&&0.5& &46& &6.76e-08     \\
  &&&&&1& &48& &5.49e-08     \\
&&&& & & & & &\\
\cline{4-10}
&&&& & &  & & &\\
 &&&&&0.005& &23& &8.74e-08      \\
&$2^{-8}$&&0.05&  &0.05& &25& &7.26e-08    \\
  &&&&&0.5& &22& &5.77e-08     \\
  &&&&&1& &23& &7.63e-08     \\
&&&& & & & & &\\
\cline{4-10}
&&&& & &  & & &\\
  &&&&&0.005& &12& &6.51e-08      \\
&&&0.5&  &0.05& &11& &8.80e-08    \\
  &&&&&0.5& &10& &7.05e-08    \\
  &&&&&1& &12& &8.53e-08     \\
&&&& & & & & &\\
\hline
\end{tabular}
\end{center}
\end{table}
\end{example}

\begin{example}\label{example:3}\rm{(Example 1 in \cite{Stadler})}
  \begin{equation*}
     \left\{ \begin{aligned}
        &\min \limits_{(y,u)\in Y\times U}^{}\ \ J(y,u)=\frac{1}{2}\|y-y_d\|_{L^2(\Omega)}^{2}+\frac{\alpha}{2}\|u\|_{L^2(\Omega)}^{2}+{\beta}\|u\|_{L^1(\Omega)} \\
        &\quad{\rm s.t.}\qquad -\Delta y=u\quad \mathrm{in}\  \Omega=(0,1)\times(0,1), \\
         &\qquad \qquad \qquad y=0\quad  \mathrm{on}\ \partial\Omega,\\
         &\qquad \qquad\qquad  u\in U_{ad}=\{v(x)|a\leq v(x)\leq b, {\rm a.e }\  \mathrm{on}\ \Omega \},
                          \end{aligned} \right.
 \end{equation*}
 where the desired state $y_d=\frac{1}{6}\sin(2\pi x)\exp(2x)\sin(2\pi y)$, the parameters $\alpha=10^{-5}$, $\beta=10^{-3}$, $a=-30$ and $b=30$. In addition, the exact solution of the problem is unknown. In this case, using a numerical solution as the reference solution is a common method. For more details, one can see \cite{HiPiUl}. In our practice implementation, we use the numerical solution computed on a grid with $h^*=2^{-10}$ as the reference solution. It should be emphasized that choosing the solution that computed on mesh $h^*=2^{-10}$ is reliable. As shown in below, when $h^*=2^{-10}$, the scale of data is 1046529.

 An example, the computed discretized optimal control $u_h$ with $h=2^{-7}$ is displayed in Figure \ref{fig:control on h=$2^{-7}$}. In Table \ref{tab:5}, we report the numerical results obtained by our proposed sGS-imABCD method for solving (\ref{eqn:discretized matrix-vector dual problem}). As a result, one can see that our proposed sGS-imABCD method is an efficient algorithm to solve problem (\ref{eqn:discretized matrix-vector dual problem}) to high accuracy. In addition, the errors of the control $u$ with respect to the solution on the finest grid ($h^*=2^{-10}$) and the results of EOC for control are also presented in Table \ref{tab:5}, which confirm the error estimate result as shown in Theorem \ref{theorem:error1}. For the sake of comparison, in Table \ref{tab:7}, we report the numerical results obtained by sGS-imABCD method for solving (\ref{eqn:discretized matrix-vector dual problem}) and iwADMM, APG methods for (\ref{equ:approx discretized matrix-vector form}). Comparing the error results from Table \ref{tab:5} and Table \ref{tab:7}, we can see that directly solving (\ref{eqn:discretized matrix-vector dual problem}) can get better error results than that from solving (\ref{eqn:discretized matrix-vector dual problem}) and (\ref{equ:approx discretized matrix-vector form}). Obviously, this conclusion show the efficiency of our dual-based approach which can avoid the additional error caused by the approximation of $L^1$-norm. Furthermore, from Table \ref{tab:5}, the numerical results in terms of iteration numbers illustrate the mesh-independent performance of our proposed sGS-imABCD method.

 In addition, in Table \ref{tab:6}, numbers of iteration steps and the relative residual errors
 of PMHSS-preconditioned GMRES method for the $\hat{p}$-subproblem on mesh $h=2^{-7}$ and $h=2^{-8}$ are presented, which shows that the PMHSS-preconditioned GMRES method is roughly independent of the mesh size $h$.

 As a result from Table \ref{tab:7}, it can be also observed that our sGS-imABCD is faster and more efficient than the iwADMM and APG methods in terms of the iteration numbers and CPU times. The numerical performance of our proposed sGS-imABCD method clearly demonstrates the importance of our method.

Finally, to show the influence of the parameters $\alpha$ and $\beta$ on our proposed sGS-imABCD method, we also test the Example \ref{example:3} with different values of $\alpha$ and $\beta$ on mesh $h=2^{-8}$. The results are presented in Table \ref{tab:8}. From the Table \ref{tab:8}, it is obviouse to see that our proposed sGS-imABCD method is independent of the parameter $\beta$. However its convergence rate depends on $\alpha$.
It also confirms the convergence results of Theorem \ref{sGS-imABCD convergence}.
\begin{figure}[H]
\begin{center}
\includegraphics[width=0.30\textwidth]{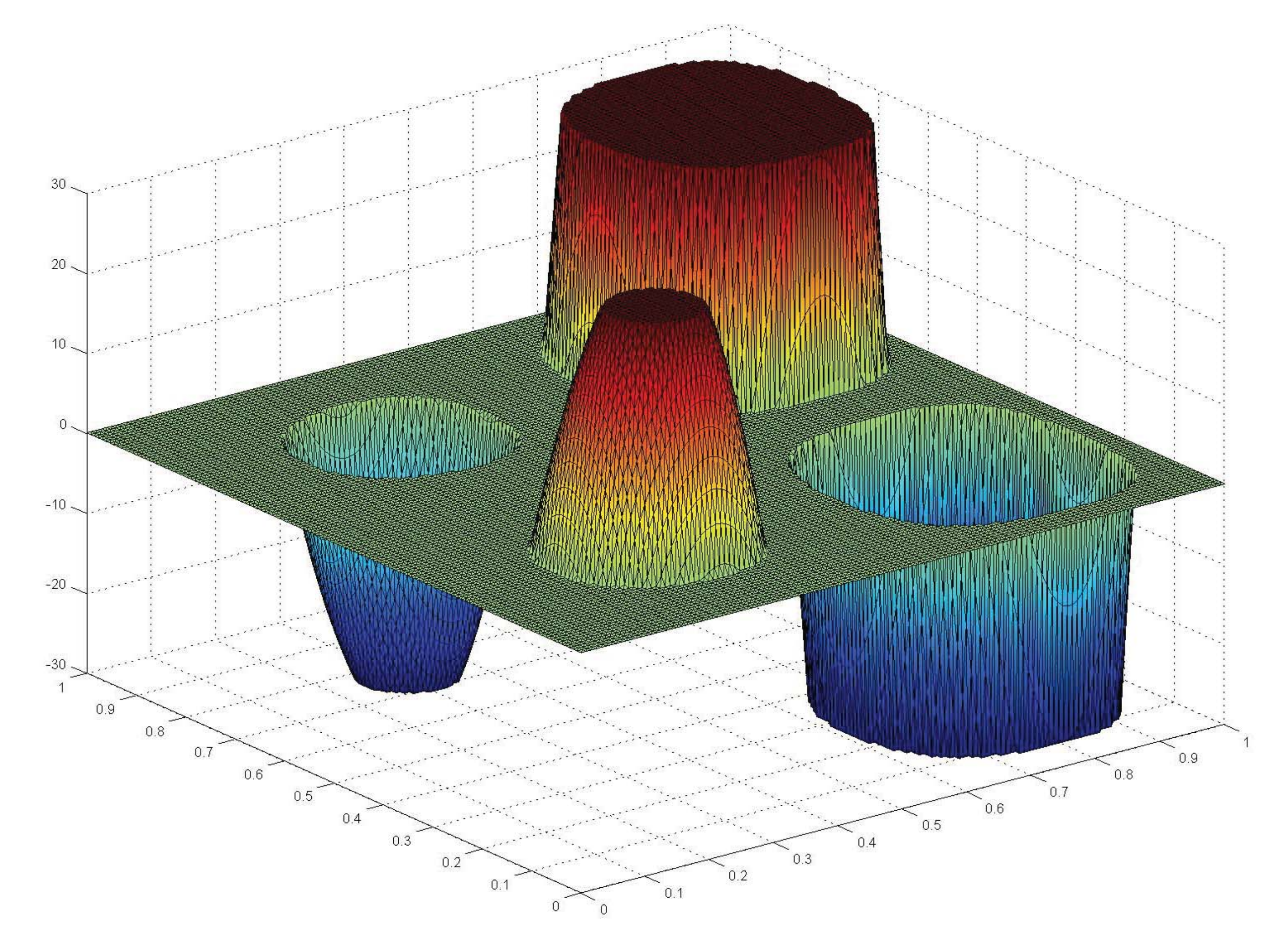}
\includegraphics[width=0.30\textwidth]{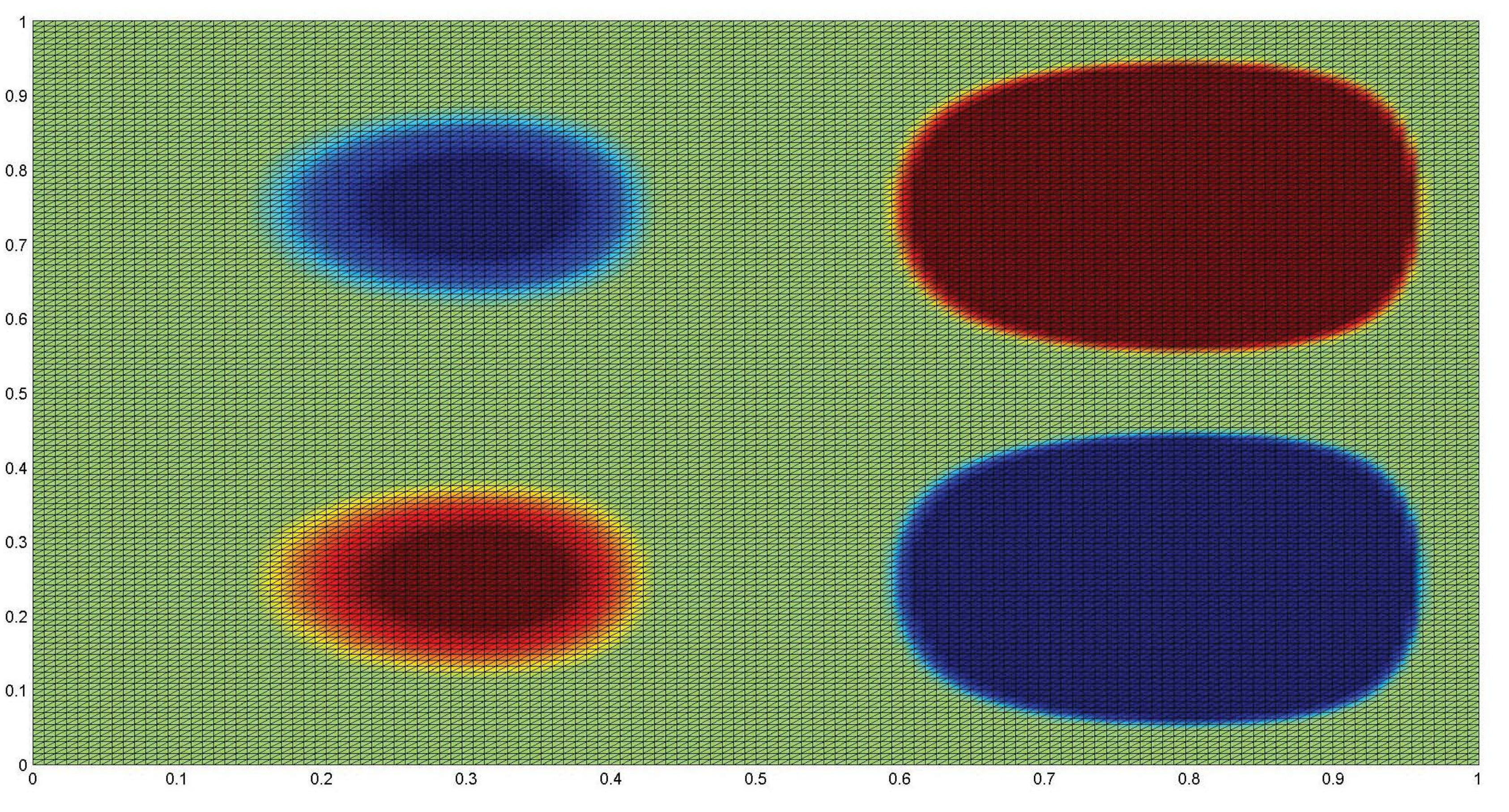}
\caption{Optimal control $u_h$ on the square, $h=2^{-7}$. Dark red and dark blue areas correspond to $u_h=\pm 30$ and green areas to $u_h=0$}\label{fig:control on h=$2^{-7}$}.
\end{center}
\end{figure}

\begin{table}[H]\footnotesize
\caption{Example \ref{example:3}: The performance of sGS-imABCD for (\ref{eqn:discretized matrix-vector dual problem}). In the table, $\#$dofs stands for the number of degrees of freedom for the control variable on each grid level.}\label{tab:5}
\begin{center}
\begin{tabular}{@{\extracolsep{\fill}}ccccccccccc}
\hline
\multirow{2}{*}{$h$}&\multirow{2}{*}{$\#$dofs} &\multirow{2}{*}{iter.sGS-imABCD} &\multirow{2}{*}{No.$\tilde{p}$-block}&\multirow{2}{*}{}& \multirow{2}{*}{residual $\eta$}&\multirow{2}{*}{} & \multirow{2}{*}{CPU time/s} & \multirow{2}{*}{$E_2$}&\multirow{2}{*}{EOC}\\
& &  &  &  & &  &  &\\
\hline
&&&&&&&\\
$2^{-3}$  &49 &37 &12   &   &8.67e-08  &   &0.64 &5.5408      &--       &    \\
&&&&&&&\\
$2^{-4}$  &225 &30 &10  &   &7.32e-08  &   &0.65 &2.4426      &1.1817       &    \\
&&&&&&&\\
$2^{-5}$  &961 &22 & 8  &   &8.38e-08  &   &0.73 &1.1504      &1.1340   &    \\
&&&&&&&\\
$2^{-6}$  &3969 &22 &7  &   &6.83e-08  &   &4.65 &0.4380    &1.2203   &    \\
&&&&&&&\\
$2^{-7}$  &16129&16 &5  &   &6.46e-08  &   &16.60 &0.1774    &1.2413  &    \\
&&&&&&&\\
$2^{-8}$  &65025&15 &3  &   &6.36e-08  &   &105.70 &0.1309   &1.0807   &    \\
&&&&&&&\\
$2^{-9}$  &261121&15 &3 &   &5.65e-08  &   &1158.62 &0.0406  & 1.1821   &    \\
&&&&&&&\\
$2^{-10}$  &1046529&16 &3&  &4.50e-08&  &24008.07 &--   &--   &    \\
&&&&&&&\\
\hline
\end{tabular}
\end{center}
\end{table}
\begin{table}[H]\footnotesize
\caption{Example \ref{example:3}: The convergence behavior of GMRES for $\hat{p}$-block subproblem. 
}\label{tab:6}
\begin{center}
\begin{tabular}{@{\extracolsep{\fill}}|c|c|c|c|}
\hline
\multirow{2}{*}{$h$}&\multirow{2}{*}{iter.sGS-imABCD}&\multirow{2}{*}{iter.GMRES of $\hat{p}$-block}&\multirow{2}{*}{Relative residual error of GMRES} \\
&&&\\
\hline

                            &\multirow{2}{*}{1}    &\multirow{2}{*}{7} &\multirow{2}{*}{1.54e-04} \\
                            &\multirow{2}{*}{2}    &\multirow{2}{*}{7} &\multirow{2}{*}{1.12e-05}\\
                            &\multirow{2}{*}{3}    &\multirow{2}{*}{8} &\multirow{2}{*}{7.25e-06}\\
                            &\multirow{2}{*}{4}    &\multirow{2}{*}{8} &\multirow{2}{*}{3.95e-06}\\
                            &\multirow{2}{*}{5}    &\multirow{2}{*}{8} &\multirow{2}{*}{3.85e-06}\\
                            &\multirow{2}{*}{6}    &\multirow{2}{*}{8} &\multirow{2}{*}{2.66e-06}\\
                            &\multirow{2}{*}{7}    &\multirow{2}{*}{8} &\multirow{2}{*}{3.33e-06}\\
  \multirow{2}{*}{$2^{-7}$} &\multirow{2}{*}{8}    &\multirow{2}{*}{8} &\multirow{2}{*}{2.60e-06}\\
                            &\multirow{2}{*}{9}   &\multirow{2}{*}{8} &\multirow{2}{*}{1.86e-06}\\
                            &\multirow{2}{*}{10}   &\multirow{2}{*}{8} &\multirow{2}{*}{1.15e-06}\\
                            &\multirow{2}{*}{11}   &\multirow{2}{*}{8} &\multirow{2}{*}{1.28e-06}\\
                            &\multirow{2}{*}{12}    &\multirow{2}{*}{7} &\multirow{2}{*}{8.68e-07}\\
                            &\multirow{2}{*}{13}   &\multirow{2}{*}{7} &\multirow{2}{*}{9.26e-07}\\
                            &\multirow{2}{*}{14}   &\multirow{2}{*}{7} &\multirow{2}{*}{5.17e-07}\\
                            &\multirow{2}{*}{15}   &\multirow{2}{*}{7} &\multirow{2}{*}{7.76e-07}\\
                            &\multirow{2}{*}{16}   &\multirow{2}{*}{7} &\multirow{2}{*}{7.39e-07}\\
                            &&&\\
\hline

                            &\multirow{2}{*}{1}    &\multirow{2}{*}{7} &\multirow{2}{*}{1.50e-04} \\
                            &\multirow{2}{*}{2}    &\multirow{2}{*}{7} &\multirow{2}{*}{1.11e-05}\\
                            &\multirow{2}{*}{3}    &\multirow{2}{*}{8} &\multirow{2}{*}{7.23e-06}\\
                            &\multirow{2}{*}{4}    &\multirow{2}{*}{8} &\multirow{2}{*}{9.61e-06}\\
\multirow{2}{*}{$2^{-8}$}   &\multirow{2}{*}{5}    &\multirow{2}{*}{9} &\multirow{2}{*}{5.56e-06}\\
                            &\multirow{2}{*}{6}    &\multirow{2}{*}{10} &\multirow{2}{*}{7.37e-07}\\
                            &\multirow{2}{*}{7}    &\multirow{2}{*}{8} &\multirow{2}{*}{3.98e-06}\\
                            &\multirow{2}{*}{8}    &\multirow{2}{*}{8} &\multirow{2}{*}{2.34e-06}\\
                            &\multirow{2}{*}{9}    &\multirow{2}{*}{8} &\multirow{2}{*}{1.96e-06}\\
                            &\multirow{2}{*}{10}    &\multirow{2}{*}{8} &\multirow{2}{*}{1.15e-06}\\
                            &\multirow{2}{*}{11}    &\multirow{2}{*}{8} &\multirow{2}{*}{1.27e-06}\\
                            &\multirow{2}{*}{12}    &\multirow{2}{*}{7} &\multirow{2}{*}{8.36e-07}\\
                            &\multirow{2}{*}{13}    &\multirow{2}{*}{7} &\multirow{2}{*}{8.16e-07}\\
                            &\multirow{2}{*}{14}    &\multirow{2}{*}{7} &\multirow{2}{*}{4.38e-07}\\
                            &\multirow{2}{*}{15}    &\multirow{2}{*}{7} &\multirow{2}{*}{7.61e-07}\\
                            &&&\\
\hline
\end{tabular}
\end{center}
\end{table}

\begin{table}[H]\footnotesize
\caption{Example \ref{example:3}: The convergence behavior of sGS-imABCD, ihADMM and APG for (\ref{equ:approx discretized matrix-vector form}). 
}\label{tab:7}
\begin{center}
\begin{tabular}{@{\extracolsep{\fill}}|c|c|c|c|c|c|c|c|}
\hline
\multirow{2}{*}{$h$}&\multirow{2}{*}{$\#$dofs}&\multirow{2}{*}{$E_2$}&\multirow{2}{*}{EOC}& \multirow{2}{*}{Index of performance} &\multirow{2}{*}{sGS-imABCD} & \multirow{2}{*}{ihADMM} & \multirow{2}{*}{APG} \\
&&&&&&&\\
\hline

                            &&&&\multirow{2}{*}{iter}            &\multirow{2}{*}{40} &\multirow{2}{*}{56} &\multirow{2}{*}{44}              \\
\multirow{2}{*}{$2^{-3}$}   &\multirow{2}{*}{49}&\multirow{2}{*}{6.6122}&\multirow{2}{*}{--}
                            &\multirow{2}{*}{residual $\eta$}
                            &\multirow{2}{*}{6.06e-08}            &\multirow{2}{*}{8.36e-08}
                            &\multirow{2}{*}{9.92e-08}             \\
                            &&&&\multirow{2}{*}{CPU time/s}          &\multirow{2}{*}{0.72} &\multirow{2}{*}{0.42} &\multirow{2}{*}{0.60}             \\
                            &&&&&&&\\
\hline
                            &&&&\multirow{2}{*}{iter}            &\multirow{2}{*}{16} &\multirow{2}{*}{55} &\multirow{2}{*}{39}             \\
\multirow{2}{*}{$2^{-4}$}   &\multirow{2}{*}{225}&\multirow{2}{*}{2.6314}&\multirow{2}{*}{1.3293}&\multirow{2}{*}{residual $\eta$}
                            &\multirow{2}{*}{9.94e-08}            &\multirow{2}{*}{9.14e-08}
                            &\multirow{2}{*}{9.74e-08}             \\
                            &&&&\multirow{2}{*}{CPU times/s}             &\multirow{2}{*}{0.48} &\multirow{2}{*}{0.62} &\multirow{2}{*}{1.03}  \\
                            &&&&&&&\\
\hline
                            &&&&\multirow{2}{*}{iter}            &\multirow{2}{*}{21} &\multirow{2}{*}{51} &\multirow{2}{*}{29}            \\
\multirow{2}{*}{$2^{-5}$}   &\multirow{2}{*}{961}&\multirow{2}{*}{1.2825}&\multirow{2}{*}{1.1831}
                            &\multirow{2}{*}{residual $\eta$}
                            &\multirow{2}{*}{5.36e-08}&\multirow{2}{*}{8.59e-08}
                            &\multirow{2}{*}{8.31e-06}             \\
                            &&&&\multirow{2}{*}{CPU time/s}          &\multirow{2}{*}{0.99} &\multirow{2}{*}{1.707} &\multirow{2}{*}{3.84}            \\
                            &&&&&&&\\
\hline
                            &&&&\multirow{2}{*}{iter}            &\multirow{2}{*}{22} &\multirow{2}{*}{46} &\multirow{2}{*}{29}              \\
\multirow{2}{*}{$2^{-6}$}   &\multirow{2}{*}{3969}&\multirow{2}{*}{0.7514}&\multirow{2}{*}{1.0458}&\multirow{2}{*}{residual $\eta$}
                            &\multirow{2}{*}{9.91e-08}&\multirow{2}{*}{6.83e-08}
                            &\multirow{2}{*}{9.38e-08}             \\
                            &&&&\multirow{2}{*}{CPU time/s}          &\multirow{2}{*}{4.95} &\multirow{2}{*}{8.34} &\multirow{2}{*}{11.94}            \\
                            &&&&&&&\\
\hline
                            &&&&\multirow{2}{*}{iter}            &\multirow{2}{*}{20} &\multirow{2}{*}{46} &\multirow{2}{*}{24}              \\
\multirow{2}{*}{$2^{-7}$}   &\multirow{2}{*}{16129}&\multirow{2}{*}{0.29304}&\multirow{2}{*}{1.1240}&\multirow{2}{*}{residual $\eta$}
                            &\multirow{2}{*}{ 9.89e-08}&\multirow{2}{*}{5.85e-08}
                            &\multirow{2}{*}{9.36e-08}             \\
                            &&&&\multirow{2}{*}{CPU time/s}          &\multirow{2}{*}{20.83} &\multirow{2}{*}{38.93} &\multirow{2}{*}{45.85}            \\
                            &&&&&&&\\
\hline
                            &&&&\multirow{2}{*}{iter}            &\multirow{2}{*}{20} &\multirow{2}{*}{48} &\multirow{2}{*}{20}              \\
\multirow{2}{*}{$2^{-8}$}   &\multirow{2}{*}{65025}&\multirow{2}{*}{0.1357}&\multirow{2}{*}{1.1213}
                            &\multirow{2}{*}{residual $\eta$}
                            &\multirow{2}{*}{4.99e-08}&\multirow{2}{*}{8.39e-08}
                            &\multirow{2}{*}{9.05e-08}             \\
                            &&&&\multirow{2}{*}{CPU time/s}          &\multirow{2}{*}{143.88} &\multirow{2}{*}{219.27} &\multirow{2}{*}{181.11}            \\
                            &&&&&&&\\
\hline
                            &&&&\multirow{2}{*}{iter}            &\multirow{2}{*}{18} &\multirow{2}{*}{50} &\multirow{2}{*}{20}              \\
\multirow{2}{*}{$2^{-9}$}   &\multirow{2}{*}{261121}&\multirow{2}{*}{0.0958}&\multirow{2}{*}{1.0181}
                            &\multirow{2}{*}{residual $\eta$}
                            &\multirow{2}{*}{9.05e-08}&\multirow{2}{*}{7.04e-08}
                            &\multirow{2}{*}{8.84e-08}             \\
                            &&&&\multirow{2}{*}{CPU time/s}          &\multirow{2}{*}{1272.25} &\multirow{2}{*}{2227.48} &\multirow{2}{*}{1959.11}            \\
                            &&&&&&&\\
\hline
\end{tabular}
\end{center}
\end{table}

\begin{table}[H]\footnotesize
\caption{Example \ref{example:3}: The performance of sGS-imABCD for (\ref{eqn:discretized matrix-vector dual problem}) with different values of $\alpha$ and $\beta$.}\label{tab:8}
\begin{center}
\begin{tabular}{@{\extracolsep{\fill}}ccccccccccc}
\hline
&\multirow{2}{*}{$h$}&\multirow{2}{*}{}&\multirow{2}{*}{$\alpha$}&\multirow{2}{*}{}&\multirow{2}{*}{$\beta$}&\multirow{2}{*}{} &\multirow{2}{*}{iter.sGS-imABCD}&\multirow{2}{*}{}&\multirow{2}{*}{residual error $\eta$ about K-K-T} \\
&&&& & &  & & &\\
\hline
&&&& & &  & & &\\
  &&&&&0.0005& &26& &8.37e-08      \\
&&&$10^{-6}$&  &0.001& &27& &8.40e-08    \\
  &&&&&0.005& &26& &9.77e-08     \\
  &&&&&0.008& &28& &2.47e-08     \\
&&&& & & & & &\\
\cline{4-10}
&&&& & &  & & &\\
  &&&&&0.0005& &13& &5.44e-08      \\
&$2^{-8}$&&$10^{-5}$&  &0.001& &15& &6.36e-08    \\
  &&&&&0.005& &14& &8.60e-08     \\
  &&&&&0.008& &13& &8.17e-08     \\
&&&& & & & & &\\
\cline{4-10}
&&&& & &  & & &\\
  &&&&&0.0005& &5& &9.84e-08      \\
&&&$10^{-4}$&  &0.001& &4& &3.71e-08    \\
  &&&&&0.005& &5& &9.23e-08     \\
  &&&&&0.008& &5& &5.22e-08     \\
&&&& & & & & &\\
\hline
\end{tabular}
\end{center}
\end{table}
\end{example}

\section{Concluding remarks}\label{sec:7}
In this paper, instead of solving the optimal control problem with $L^1$ control cost, we directly solve the dual problem which is an unconstrained multi-block minimization problem. By taking advantage of the structure of dual problem, and combining the inexact majorized ABCD (imABCD) method and the recent advances in the inexact symmetric Gauss-Seidel (sGS) technique, we introduce the sGS-imABCD method to solve the dual problem. Its efficiency is confirmed by both the theory and numerical results. As it is mentioned, the iwADMM and APG methods could be employed to solve (\ref{equ:approx discretized matrix-vector form}), the approximative discretization of the primal problem (\ref{eqn:orginal problems}). For the sake of comparison, we also use our method to solve its dual ({\ref{eqn:approx discretized matrix-vector dual problem}}). As shown in the numerical results, directly solving the dual problem (\ref{eqn:discretized matrix-vector dual problem}) could get better error results than that from solving (\ref{equ:approx discretized matrix-vector form}) and ({\ref{eqn:approx discretized matrix-vector dual problem}}). It should be stressed that the better error results are due to the fact that solving (\ref{eqn:discretized matrix-vector dual problem}) can avoid the approximation of the discrete $L^1$-norm. More importantly, numerical experiments show that our proposed method for solving ({\ref{eqn:approx discretized matrix-vector dual problem}}) outperforms the ihADMM and APG methods for solving (\ref{equ:approx discretized matrix-vector form}).

\section*{Acknowledgments}
The authors would like to thank Prof. Defeng Sun and Prof. Kim-Chuan Toh at National University of Singapore for their valuable suggestions that led to improvement in this paper and also would like to thank Prof. Long Chen for the FEM package iFEM \cite{Chen} in Matlab.

\end{document}